\documentclass[11pt]{article}
\usepackage{amssymb,amsmath,latexsym, epsfig}

\begin{document}

\begin{doublespace}

\newtheorem{thm}{Theorem}[section]
\newtheorem{lemma}[thm]{Lemma}
\newtheorem{defn}{Definition}[section]
\newtheorem{prop}[thm]{Proposition}
\newtheorem{corollary}[thm]{Corollary}
\newtheorem{remark}[thm]{Remark}
\newtheorem{example}[thm]{Example}
\numberwithin{equation}{section}

\def\ee{\varepsilon}
\def\qed{{\hfill $\Box$ \bigskip}}
\def\MM{{\cal M}}
\def\BB{{\cal B}}
\def\LL{{\cal L}}
\def\FF{{\cal F}}
\def\GG{{\cal G}}
\def\EE{{\cal E}}
\def\QQ{{\cal Q}}

\def\R{{\mathbb R}}
\def\L{{\bf L}}
\def\E{{\mathbb E}}
\def\F{{\bf F}}
\def\P{{\mathbb P}}
\def\N{{\mathbb N}}
\def\eps{\varepsilon}
\def\wh{\widehat}
\def\pf{\noindent{\bf Proof.} }

\title{\Large \bf On the Estimates of the Density of the Purely Discontinuous Girsanov Transforms of $\alpha$-Stable-like Processes          }
\author{ Chunlin Wang\\
Department of Mathematics\\
University of Illinois\\
Urbana, IL 61801, USA\\
 Email: cwang13@uiuc.edu\\   }
\date{}
\maketitle

\begin{abstract}
In this paper, we study the purely discontinuous Girsanov
transforms which were discussed in Chen and Song \cite{CS2}
 and Song \cite{S3}. We show that the transition density of any
 purely discontinuous Girsanov transform of a $\alpha$-stable-like
process, which can be nonsymmetric, is comparable to the
transition density of the $\alpha$-stable-like process.
\end{abstract}



 \vspace{.1truein}
\noindent{\bf Running Title:} Density estimates of the purely
discontinuous Girsanov transformed processes

\pagebreak
\section{Introduction}
Suppose that $(X,\mathbb{P}_x)$ is a Brownian motion in $\R^d$ and
${\cal M}_t$ is the $\sigma$-field generated by $\{X_s,s\leq t\}$.
Define
$$L_t=\exp\left(\int_{0}^{t}b(X_s)dX_s-\frac{1}{2}\int_{0}^{t}|b(X_s)|^2ds\right),$$
where $b$ is an $\R^d$-valued function on $\R^d$. When $b$
satisfies certain natural conditions, $L_t$ is a nonnegative local
martingale under each $\mathbb{P}_x$ and thus is a supermatingale
multiplicative additive functional of $X$. $L_t$ defines a family
of probability measure $\{\tilde{\mathbb{P}}_x:x\in \R^d\}$ on
${\cal M}_{\infty}$ by $d\tilde{\mathbb{P}}_x=L_td\mathbb{P}_x$ on
${\cal M}_t$. Under this family of measures
$\{\tilde{\mathbb{P}}_x:x\in \R^d\}$, $X_t$ is denoted by
$\tilde{X}_t$. The process $\tilde{X}$ is called a Girsanov
transform of $X$ and the transition density of $\tilde{X}$ has
Gaussian upper and lower estimates for certain natural function
$b$.

 The main purpose of this paper is to study upper and lower
 estimates on the transition densities of Girsanov transforms of
 $\alpha$-stable-like processes (see definition 1.1).

Let's first recall the purely discontinuous Girsanov transforms
studied in \cite{CS2}. For convenience, we consider the state
space $\R^d$ instead of a Lusin space $E$ in \cite{CS2}. Let
${\cal{B}}(\R^d)$ be the Borel $\sigma$-field and $m$ a
$\sigma$-finite measure on ${\cal{B}}(\R^d)$ . Let
$X=(X_t,\mathbb{P}_t)$ be a $m$-Hunt process on $\R^d$ and has a
L\'{e}vy system $(N,H)$. That is for any nonnegative function $f$
on $\R^d\times\R^d$ vanishing on the diagonal
$$\mathbb{E}_x\left(\sum_{s \le
t}f(X_{s-},X_{s})\right)=\mathbb{E}_x\left(\int_{0}^{t}\int_{R^d}f(X_s,y)N(X_s,dy)dH_s\right),$$
for every $x \in \R^d$ and $t>0$.

 Let $F$ be a bounded function on
$\R^d\times\R^d$ that vanishes on the diagonal. We say
$F\in{\bf{J}}(X)$ if
$$\lim_{t\downarrow 0}\sup_{x\in
\R^d}\mathbb{E}_x\left(\int_{0}^{t}\int_{R^d}|F|(X_s,y)N(X_s,dy)dH_s\right)=0,$$
see \cite{CS}. If $F\in {\bf{J}}(X)$ also satisfies the condition
$\inf_{x,y\in E}F(x,y)>-1$, then the process
$$t\mapsto
\sum_{0<s<t}F(X_{s-},X_s)-\int_{0}^{t}\int_{R^d}F(X_s,y)N(X_s,dy)dH_s$$
is a martingale and its Doleans-Dade exponential is (see Theorem
9.39 of \cite{HWY})
\begin{eqnarray*}
      L_t&=&\exp\left(-\int_{0}^{t}\int_{R^d}F(X_s,y)N(X_s,dy)dH_s\right)\prod_{s\leq t}(1+F(X_{s-},X_s))\\
         &=&\exp\left(\sum_{s\leq t}\ln(1+F(X_{s-},X_s))-\int_{0}^{t}\int_{R^d}F(X_s,y)N(X_s,dy)dH_s\right).\\
\end{eqnarray*}
$L_t$ is a nonnegative local martingale and thus a supermartingale
multiplicative functional of $X$. $L_t$ defines a family of
probability of measures $\{\mathbb{P}_{x}, x\in \R^d\}$ on ${\cal
M}_{\infty}$ by $d\tilde{\mathbb{P}}_{x}=L_t d \mathbb{P}_{x}$ on
${\cal M}_{t}$, see section 62 of \cite{S}. Let
$\tilde{X}=(\tilde{X}_t, \cal{M},$$\,{\cal M}_t,
\tilde{\mathbb{P}}_x, x \in {\R^d} )$ denote this new process.
 Under this family of measures $\tilde{\mathbb{P}}_{x}$, $X_t$ is denoted by
$\tilde{X}_t$. The process $\tilde{X}_t$ is called a purely
discontinuous Girsanov transform of $X$. Under further assumption
that $F\in {\bf{A}}_2(X)$ (see \cite{CS2} for definition) and with
a Lusin space $E$ instead of $\R^d$, it was shown in \cite{CS2}
that the Green function of a purely discontinuous Girsanov
transform of $\tilde{X}$ is comparable to that of $X$, especially
when $X$ is a symmetric stable process. Under the assumption that
$F\in {\bf{I}}_2(X)$ (see \cite{S3} for definition) and with a
Lusin space $E$ instead of $\R^d$, it was shown in \cite{S3} that
when $X$ is a symmetric stable process, the transition density of
a purely discontinuous Girsanov transform of $\tilde{X}$ is
comparable to that of $X$.

One of the main tools used in \cite{CS2} and \cite{S3} is the
Dirichlet form. The Dirichlet form was also used in \cite{S2} to
study the two-sided estimates on the density of the Feynman-Kac
semigroups of symmetric $\alpha$-stable-like processes (see
definition 1.1). The Dirichlet form works well when the processes
are symmetric. It encounters difficulty when the processes are not
symmetric.

Let's introduce $\alpha$-stable-like processes in the following.

Suppose the process $X=(X_t,\mathbb{P}_t)$ given above has a
L\'{e}vy system $(N,H)$ given by $H_{t}=t$ and
$$ N(x,dy)=2C(x,y){|x-y|}^{-(d+\alpha)}m(dy),$$
where $m$ a $\sigma$-finite measure on ${\cal{B}}(\R^d)$ given by
$m(dx)=M(x)dx$ with $M(x)$ nonnegative and bounded.

\begin{defn}
We say that $X$ is an $\alpha$-stable-like process if $C(x,y)$ is
bounded.
\end{defn}

In this paper we assume that $X$ admits a transition density
$p(t,x,y)$ with respect to $m$ and $p(t,x,y)$ is jointly
continuously on $(0,\infty)\times\R^d\times\R^d$ and satisfies the
condition
\begin{equation}
 M_{1}t^{-\frac{d}{\alpha}}\left(1 \wedge
\frac{t^{\frac{1}{\alpha}}}{|x-y|}\right)^{d+\alpha}\leq
p(t,x,y)\leq M_{2}t^{-\frac{d}{\alpha}}\left(1 \wedge
\frac{t^{\frac{1}{\alpha}}}{|x-y|}\right)^{d+\alpha},\,\,\forall
(t,x,z)\in (0,\infty)\times\R^d\times\R^d,
\end{equation}
where $M_{1}$ and $M_{2}$ are positive constants.

Here we do not assume that $X$ is symmetric. When $X$ is
symmetric, it is called a symmetric $\alpha$-stable-like process,
which was introduced in \cite{CT}, where a symmetric Hunt process
is associated with a regular Dirichlet form and thus Dirichlet
form method can be applied.

We list some examples for $\alpha$-stable-like processes. One
dimensional $\alpha$-stable processes with their L\'{e}vy measure
$\nu$ concentrated neither on $(0,\infty)$ nor on $(-\infty,0)$
are $\alpha$-stable-like processes, as their $C(x,y)$ are bounded
and densities $p(t,x,y)$ satisfies (1.1), see  Definition 14.16
and Remark 14.18 in \cite{SK}. For higher dimensions, \cite{VZ}
gave a large class of nonsymmetric strictly $\alpha$-stable
processes which satisfy (1.1) with $C(x,y)$ bounded. So this class
also belongs to the class of $\alpha$-stable-like processes.

In \cite{RW}, we developed a method, which was based on an idea in
\cite{BM} for Brownian motions, some results on discontinuous
functionals and some techniques for integral estimates, to obtain
the estimates on the density of the Feynman-Kac semigroups of
$\alpha$-stable-like processes where neither $F$ nor $X$ is
symmetric. In this paper, by improving this method so that it also
works for the purely discontinuous Gisarnov transform, we will
show that the transition density of the purely discontinuous
Girsanov transform of $\tilde{X}$ is comparable to that of $X$
with neither $F$ nor $X$ being symmetric.

The content of this paper is organized as follows. In section 2,
we present preliminary results on additive functionals. In section
3, we establish the two-sided estimates on the transition density
of the purely discontinuous Girsanov transformations under certain
assumptions of $F(x,y)$.

\section{Preliminary Results on Additive Functionals}

For any $l>0$, let $A^{(l)}_t=B^{(l)}_t-D^{(l)}_t$, where
$$B^{(l)}_t=\sum_{s \le
t}\ln(1+F(X_{s-},X_{s}))1_{\{|X_{s-}-X_{s}|>1/l\}},$$ and
$$D^{(l)}_t=\int_{0}^{t}\int_{E}F(X_s,y)1_{\{|X_s-y|>1/l\}}N(X_s,dy)dH_s, \,\, l=1,2,\cdots.$$
$B^{(l)}_t$ is a pure discontinuous functional and $D^{(l)}_t$ is
a continuous functional. For convenience of notation, we omit the
index $l$, i.e. denote $A^{(l)}_t,\, B^{(l)}_t$ and $D^{(l)}_t$ by
$A_t, \,B_t$ and $D_t$ respectively. We have the following
formulae for $A^{n}_t$.

\begin{thm}
\begin{eqnarray*}
            A^n_{t}&=&C^1_{n}\int_{0}^{t}A_{s}^{n-1}\,dA_{s}-C^2_{n}\int_{0}^{t}A_{s}^{n-2}(\ln(1+F(X_{s-},X_{s})))1_{\{|X_{s-}-X_{s}|>1/l\}}\,dA_{s}+\cdots\\
            &&+(-1)^{i-1}C^{i}_n \int_{0}^{t}A_{s}^{n-i}(\ln(1+F(X_{s-},X_{s})))^{i-1}1_{\{|X_{s-}-X_{s}|>1/l\}}
            \,dA_{s}+\cdots\\
            &&+(-1)^{n-1}C^n_n\int_{0}^{t}(\ln(1+F(X_{s-},X_{s})))^{n-1}1_{\{|X_{s-}-X_{s}|>1/l\}} \,dA_{s}\,\\
\end{eqnarray*}
and
\begin{eqnarray*}
          A^n_{t}&=&C^1_{n}\int_{0}^{t}(A_{t}-A_{s})^{n-1}\,dA_{s}+C^2_{n}\int_{0}^{t}(A_{t}-A_{s})^{n-2}(\ln(1+F(X_{s-},X_{s})))1_{\{|X_{s-}-X_{s}|>1/l\}} \,dA_{s}\\
                 &&+\cdots+C^{i}_n \int_{0}^{t}(A_{t}-A_{s})^{n-i}(\ln(1+F(X_{s-},X_{s})))^{i-1}1_{\{|X_{s-}-X_{s}|>1/l\}}\,dA_{s}+\cdots\\
                 &&+C^n_n\int_{0}^{t}(\ln(1+F(X_{s-},X_{s})))^{n-1}1_{\{|X_{s-}-X_{s}|>1/l\}}\,dA_{s}.\\
\end{eqnarray*}
\end{thm}
\pf Note that
$$
A_{s}-A_{s-}=B_{s}-B_{s-}=\ln(1+F(X_{s-},X_{s}))1_{\{|X_{s-}-X_{s}|>1/l\}}.
$$
 We use induction to show these two formulae for $A^n_t$. It is
clear that they hold for $n=2$. Suppose they hold for $n \leq
m-1$, we show they hold for $n=m$.

It follows from the integration by parts formula,
$$A^m_t=\int_{0}^{t}A_{s-}\,dA^{m-1}_{s}+\int_{0}^{t}A_{s}^{m-1}\,dA_{s},$$
where
\begin{eqnarray*}
                &&\int_{0}^{t}A_{s-}\,dA^{m-1}_{s}\\
                &=&\int_{0}^{t}(A_s-\ln(1+F(X_{s-},X_{s}))1_{\{|X_{s-}-X_{s}|>1/l\}})\,dA^{m-1}_{s}\\
                &=&\int_{0}^{t}A_s\,dA^{m-1}_{s}-\int_{0}^{t}\ln(1+F(X_{s-},X_{s}))1_{\{|X_{s-}-X_{s}|>1/l\}}\,dA^{m-1}_{s}\\
                &=&\int_{0}^{t}A_s(\sum_{i=1}^{m-1}(-1)^{i-1}C_{m-1}^{i}A^{m-1-i}_{s}(\ln(1+F(X_{s-},X_{s}))1_{\{|X_{s-}-X_{s}|>1/l\}})^{i-1})\,dA_{s}\\
                &&-\int_{0}^{t}F(X_{s-},X_{s})(\sum_{j=1}^{m-1}(-1)^{j-1}C_{m-1}^{j}A^{m-1-j}_{s}(\ln(1+F(X_{s-},X_{s}))1_{\{|X_{s-}-X_{s}|>1/l\}})^{j-1}\,\\
                && \cdot dA_{s}\,\, (\textrm{ by the first formula for } A^n_t \textrm{ when } n=m-1 \textrm{ } ) \\
                &=&\sum_{i=1}^{m-1}(-1)^{i-1}C_{m-1}^{i}\int_{0}^{t}A^{m-i}_{s}(\ln(1+F(X_{s-},X_{s}))1_{\{|X_{s-}-X_{s}|>1/l\}})^{i-1}\,dA_{s}\\
                &&-\sum_{j=1}^{m-1}(-1)^{j-1}C_{m-1}^{j}\int_{0}^{t}A^{m-1-j}_{s}(\ln(1+F(X_{s-},X_{s}))1_{\{|X_{s-}-X_{s}|>1/l\}})^{j-1}\,dA_{s}\\
\end{eqnarray*}
\begin{eqnarray*}
                &=&\sum_{i=1}^{m-1}(-1)^{i-1}C_{m-1}^{i}\int_{0}^{t}A^{m-i}_{s}(\ln(1+F(X_{s-},X_{s}))1_{\{|X_{s-}-X_{s}|>1/l\}})^{i-1}\,dA_{s}\\
                &&-\sum_{i=2}^{m}(-1)^{i-2}C_{m-1}^{i-1}\int_{0}^{t}A^{m-i}_{s}(\ln(1+F(X_{s-},X_{s}))1_{\{|X_{s-}-X_{s}|>1/l\}})^{i-1}\,dA_{s}\\
                &&\textrm{                      } (\textrm{ let }j=i-1 )\\
                &=&\sum_{i=2}^{m-1}(-1)^{i-1}(C_{m-1}^{i}+C_{m-1}^{i-1})\int_{0}^{t}A^{m-i}_{s}(\ln(1+F(X_{s-},X_{s}))1_{\{|X_{s-}-X_{s}|>1/l\}})^{i-1}\,dA_{s}\\
                &&+C_{m-1}^{1}\int_{0}^{t}A_{s}^{m-1}\,dA_{s}-(-1)^{m-2}\int_{0}^{t}(\ln(1+F(X_{s-},X_{s}))1_{\{|X_{s-}-X_{s}|>1/l\}})^{m-1}\,dA_{s}\\
                &=&\sum_{i=2}^{m-1}(-1)^{i-1}C_{m}^{i}\int_{0}^{t}A^{m-i}_{s}(\ln(1+F(X_{s-},X_{s}))1_{\{|X_{s-}-X_{s}|>1/l\}})^{i-1}\,dA_{s}\\
                &&+C_{m-1}^{1}\int_{0}^{t}A_{s}^{m-1}\,dA_{s}-(-1)^{m}\int_{0}^{t}(\ln(1+F(X_{s-},X_{s}))1_{\{|X_{s-}-X_{s}|>1/l\}})^{m-1}\,dA_{s}\\
                &&\textrm{   } (\textrm{ by }C_{m-1}^{i}+C_{m-1}^{i-1}=C_{m}^{i}
                ).
\end{eqnarray*}
Thus
\begin{eqnarray*}
A^m_t&=&\int_{0}^{t}A_{s-}\,dA^{m-1}_{s}+\int_{0}^{t}A_{s}^{m-1}\,dA_{s}\\
      &=&\sum_{i=1}^{m}(-1)^{i-1}C_{m}^{i}\int_{0}^{t}A^{m-i}_{s}(\ln(1+F(X_{s-},X_{s}))1_{\{|X_{s-}-X_{s}|>1/l\}})^{i-1}\,dA_{s},
\end{eqnarray*}
i.e. the first formula for $A^n_t$ holds for n=m.

Now we go to the second formula for $A^n_t$, for $n=m$.
\begin{eqnarray*}
     &&C^1_{m}\int_{0}^{t}(A_{t}-A_{s})^{m-1}\,dA_{s}\\
     &=&C^1_{m}\int_{0}^{t}\sum_{i=0}^{m-1}C_{m-1}^{i}A^{i}_{t}(-1)^{m-1-i}A^{m-1-i}_{s}\,dA_{s}\\
     &=&\sum_{i=0}^{m-1}(-1)^{m-1-i}C^1_{m}C_{m-1}^{i}A^{i}_{t}\int_{0}^{t}A^{m-1-i}_{s}\,dA_{s}\\
      &=&\sum_{i=0}^{m-1}(-1)^{m-1-i}C_{m}^{i}(m-i)A^{i}_{t}\int_{0}^{t}A^{m-1-i}_{s}\,dA_{s}\\
      &&\textrm{        } (\textrm{ by }C_{m}^{1}C_{m-1}^{i}=C_{m}^{i}(m-i)\textrm{ } )\\
      &=&\sum_{i=0}^{m-1}(-1)^{m-1-i}C_{m}^{i}A^{i}_{t}((m-i)\int_{0}^{t}A^{m-1-i}_{s}\,dA_{s})\\
\end{eqnarray*}
\begin{eqnarray*}
      &=&\sum_{i=0}^{m-1}(-1)^{m-1-i}C_{m}^{i}A^{i}_{t}(A^{m-i}_{t}+
      \int_{0}^{t}\sum_{k=2}^{m-i}(-1)^{k}C_{m-i}^{k}A^{m-i-k}_{s}(\ln(1+F(X_{s-},X_{s}))\\
      &&\cdot 1_{\{|X_{s-}-X_{s}|>1/l\}})^{k-1}dA_{s})\\
      &&\textrm {        }(\textrm{ by the first formula of }A^n_t \textrm{ for } n=m-i \textrm{ } )\\
      &=&\sum_{i=0}^{m-1}(-1)^{m-1-i}C_{m}^{i}A^{i}_{t}A^{m-i}_{t}+\int_{0}^{t}\sum_{i=0}^{m-1}(-1)^{m-1-i}\sum_{k=2}^{m-i}(-1)^{k}C_{m}^{i}C_{m-i}^{k}A^{i}_{t}A^{m-i-k}_{s}\\
      &&\cdot (\ln(1+F(X_{s-},X_{s}))1_{\{|X_{s-}-X_{s}|>1/l\}})^{k-1}\,dA_{s},\\
\end{eqnarray*}
where
$$\sum_{i=0}^{m-1}(-1)^{m-1-i}C_{m}^{i}A^{i}_{t}A^{m-i}_{t}=(\sum_{i=0}^{m-1}(-1)^{m-1-i}C_{m}^{i})A^{m}_{t}=(1)A^{m}_{t},$$
and
\begin{eqnarray*}
              &&\int_{0}^{t}\sum_{i=0}^{m-1}(-1)^{m-1-i}\sum_{k=2}^{m-i}(-1)^{k}C_{m}^{i}C_{m-i}^{k}A^{i}_{t}A^{m-i-k}_{s}(\ln(1+F(X_{s-},X_{s}))1_{\{|X_{s-}-X_{s}|>1/l\}})^{k-1}\,dA_{s}\\
              &=&\int_{0}^{t}\sum_{k=2}^{m}\sum_{i=0}^{m-k}(-1)^{m-k-i-1}C_{m}^{k}C_{m-k}^{i}A^{i}_{t}A^{m-k-i}_{s}(\ln(1+F(X_{s-},X_{s}))1_{\{|X_{s-}-X_{s}|>1/l\}})^{k-1}\,dA_{s}\\
              &&\textrm{      }(\textrm{ by }C_{m}^{i}C_{m-i}^{k}=C_{m}^{k}C_{m-k}^{i} \textrm{ and } (-1)^{m-1-i+k}=(-1)^{m-k-i-1}\textrm{ } )\\
              &=&\sum_{k=2}^{m}C_{m}^{k}(-1)^{-1}\int_{0}^{t}(A_{t}-A_{s})^{m-k}(\ln(1+F(X_{s-},X_{s}))1_{\{|X_{s-}-X_{s}|>1/l\}})^{k-1}\,dA_{s}\\
              &=&-\sum_{k=2}^{m}C_{m}^{k}\int_{0}^{t}(A_{t}-A_{s})^{m-k}(\ln(1+F(X_{s-},X_{s}))1_{\{|X_{s-}-X_{s}|>1/l\}})^{k-1}\,dA_{s},\\
\end{eqnarray*}
therefore
\begin{eqnarray*}
  &&C^1_{m}\int_{0}^{t}(A_{t}-A_{s})^{m-1}\,dA_{s}\\
  &=&A^{m}_{t}-\sum_{k=2}^{m}C_{m}^{k}\int_{0}^{t}(A_{t}-A_{s})^{m-k}(\ln(1+F(X_{s-},X_{s}))1_{\{|X_{s-}-X_{s}|>1/l\}})^{k-1}\,dA_{s},\\
\end{eqnarray*}
i.e. the second formula for $A^n_t$ holds for $n=m$. \qed

\section{Transition Density of Pure Jump Girsanov Transformation}

We need the following lemmas to prove the main results.
\begin{lemma}
For any two positive constants $K<1$ and $L$, there exist two
constants $C_{01}(K,L)$ and $C_{02}(K,L)$ which depend on $K$ and
$L$, such that
\begin{eqnarray*}
      &&L^{-1}K^{n-1}+\frac{K^{n-2}}{2!}+\frac{LK^{n-3}}{3!}+\cdots+\frac{L^{i-2}K^{n-i}}{i!}+\cdots+\frac{L^{n-2}}{n!}\leq C_{01}(K,L)K^n,\\
      &&LK^{n-1}+\frac{L^2K^{n-2}}{2!}+\frac{L^3K^{n-3}}{3!}+\cdots+\frac{L^{i}K^{n-i}}{i!}+\cdots+\frac{L^{n}}{n!}\leq C_{02}(K,L)K^n,\,\, \textrm{ for any }n \ge 1.\\
\end{eqnarray*}

\end{lemma}
\pf Use the fact that there exists $ i_{0}\ge0$, such that when
$j\ge i_{0}$,
$$\frac{L^{j}}{j!}\leq \left(\frac{K}{2}\right)^{j}.$$

\qed
\begin{lemma}
For any two positive constants $K<1$ and $L$, there exists a
constant $C_{03}(K,L)$ which depends on $K$ and $L$, such that
$$\sum_{i=2}^{m}\sum_{j=2}^{m-i}C^{i}_{m}C^{j}_{m-i}(m-i-j)!L^{i-2}L^{j-2}K^{m-i-j}\leq
      C_{03}(K,L)m!K^m, \,\, \textrm{ for all } m \ge 2. $$

\end{lemma}
\pf
\begin{eqnarray*}
      &&\sum_{i=2}^{m}\sum_{j=2}^{m-i}C^{i}_{m}C^{j}_{m-i}(m-i-j)!L^{i-2}L^{j-2}K^{m-i-j}\\
      &&=\sum_{i=2}^{m}\sum_{j=2}^{m-i}\frac{m!}{i!j!}L^{i-2}L^{j-2}K^{m-i-j}\\
      &&=m!\sum_{i=2}^{m}\frac{L^{i-2}}{i!}\left(\sum_{j=2}^{m-i}\frac{L^{j-2}}{j!}K^{m-i-j}\right)\\
      &&\leq m!\sum_{i=2}^{m}\frac{L^{i-2}}{i!}\tilde{C}_{01}(L,K)K^{m-i}\\
      &&\textrm{ }\textrm{ } (\textrm{ by the same argument of lemma 3.1, where }\tilde{C}_{01}(L,K)\textrm{ is a positive constant }) \\
      &&\leq C_{03}(L,K)m!K^m\\
      &&\textrm{ }\textrm{ } (\textrm{ by the same argument of lemma 3.1, where }\tilde{C}_{03}(L,K) \textrm{ is a positive constant }). \\
\end{eqnarray*}
\qed

From now on we define $q_{0}(t,x,y)=p(t,x,y)$ where $p(t,x,y)$ is
the transition density of $\alpha$-stable-like process $X$ and
satisfies (1.1).

By the second formula for $A^n_t$, we have for any $g$ bounded
measurable
\begin{eqnarray*}
       &&\mathbb{E}_x[A^n_{t}g(X_t)]\\
       &&=\sum_{i=1}^nC^{i}_n\mathbb{E}_x[\int_{0}^{t}(A_{t}-A_{s})^{n-i}g(X_t)(\ln(1+F(X_{s-},X_{s})))^{i-1}1_{\{|X_{s-}-X_{s}|>1/l\}}\,dA_{s}]\\
       &&=\sum_{i=1}^nC^{i}_n\mathbb{E}_x[\int_{0}^{t}(A_{t}-A_{s})^{n-i}g(X_t)(\ln(1+F(X_{s-},X_{s})))^{i-1}1_{\{|X_{s-}-X_{s}|>1/l\}}\,dB_{s}]\\
       &&\textrm{ }\textrm{ }\,\,-\sum_{i=1}^nC^{i}_n\mathbb{E}_x[\int_{0}^{t}(A_{t}-A_{s})^{n-i}g(X_t)(\ln(1+F(X_{s-},X_{s})))^{i-1}1_{\{|X_{s-}-X_{s}|>1/l\}}\,dD_{s}]\\
       &&=\sum_{i=1}^nC^{i}_n\mathbb{E}_x[\int_{0}^{t}\mathbb{E}_{X_s}\left(A_{t-s}^{n-i}g(X_{t-s})\right)\,d(\sum_{r
\le
                                  s}(\ln(1+F(X_{r-},X_{r})))^{i}1_{\{|X_{r-}-X_{r}|>1/l\}}]\\
       &&\textrm{  }\textrm{ }\,-C^{1}_{n}\mathbb{E}_x[\int_{0}^{t}\mathbb{E}_{X_s}\left(A_{t-s}^{n-1}g(X_{t-s})\right)\,dD_s]\\
       &&=\sum_{i=1}^nC^{i}_n\mathbb{E}_x[\int_{0}^{t}\int_{\R^d}\frac{2C(X_s,y)(\ln(1+F(X_{s},y)))^{i-1}1_{\{|X_{s}-y|>1/l\}}}{|X_s-y|^{d+\alpha}}\mathbb{E}_{y}\left(A_{t-s}^{n-i}g(X_{t-s})\right)\,m(dy)ds]\\
       &&\textrm{  }\textrm{ }\,-C^{1}_{n}\mathbb{E}_x[\int_{0}^{t}\int_{\R^d}\frac{2C(X_s,y)F(X_{s},y)1_{\{|X_{s}-y|>1/l\}}}{|X_s-y|^{d+\alpha}}\mathbb{E}_{y}\left(A_{t-s}^{n-1}g(X_{t-s})\right)\,m(dy)ds]\\
       &&=C^{1}_{n}\mathbb{E}_x[\int_{0}^{t}\int_{\R^d}\frac{2C(X_s,y)(\ln(1+F(X_{s},y))-F(X_{s},y))1_{\{|X_{s}-y|>1/l\}}}{|X_s-y|^{d+\alpha}}\mathbb{E}_{y}\left(A_{t-s}^{n-1}g(X_{t-s})\right)\,m(dy)\\
       &&\textrm{  }\textrm{ }\cdot ds]+\sum_{i=2}^nC^{i}_n\mathbb{E}_x[\int_{0}^{t}\int_{\R^d}\frac{2C(X_s,y)(\ln(1+F(X_{s},y)))^{i-1}1_{\{|X_{s}-y|>1/l\}}}{|X_s-y|^{d+\alpha}}\mathbb{E}_{y}\left(A_{t-s}^{n-i}g(X_{t-s})\right)\,m(dy)\\
       &&\textrm{  }\textrm{ }\cdot ds].\\
 \end{eqnarray*}
We define $q^{(l)}_n(t,x,z)$ as follows,
\begin{eqnarray*}
       &&q^{(l)}_n(t,x,z)\\
       &&=C^{1}_n\int_{0}^{t}\int_{\R^d}p(s,x,w)m(dw)\int_{\R^d}\frac{2C(w,y)(\ln(1+F(w,y))-F(w,y))1_{\{|w-y|>1/l\}}}{|w-y|^{d+\alpha}}\\
       &&\textrm{ }\textrm{ }\cdot q_{n-1}(t-s,y,z)m(dy)ds +\sum_{i=2}^nC^{i}_n\int_{0}^{t}\int_{\R^d}p(s,x,w)m(dw)\int_{\R^d}\frac{2C(w,y)(\ln(1+F(w,y)))^{i}}{|w-y|^{d+\alpha}} \\
       &&\textrm{ }\textrm{ }\cdot 1_{\{|w-y|>1/l\}}q_{n-i}(t-s,y,z)\,m(dy)ds.\\
\end{eqnarray*}
For convenience, we denote $q^{(l)}_n(t,x,z)$ by $q_n(t,x,z).$

Then by induction, we can show that for any $n \geq 1$,
$$\int_{\R^d}q_n(t,x,z)g(z)\,m(dz)=\mathbb{E}_x[A^n_{t}g(X_t)],$$
and
 \begin{eqnarray*}
       &&\mathbb{E}_x[A^n_{t}g(X_t)]\\
       &&=C^{1}_n\mathbb{E}_x[\int_{0}^{t}\int_{\R^d}\frac{2C(X_s,y)(\ln(1+F(X_{s},y))-F(X_{s},y))1_{\{|X_{s}-y|>1/l\}}}{|X_s-y|^{d+\alpha}}\int_{\R^d}q_{n-1}(t-s,y,z)g(z)\\
       &&\textrm{ }\textrm{ }\cdot m(dz)\,m(dy)ds]+\sum_{i=2}^nC^{i}_n\mathbb{E}_x[\int_{0}^{t}\int_{\R^d}\frac{2C(X_s,y)(\ln(1+F(X_{s},y)))^{i}1_{\{|X_{s}-y|>1/l\}}}{|X_s-y|^{d+\alpha}}\\
       &&\textrm{ }\textrm{ }\cdot \int_{\R^d}q_{n-i}(t-s,y,z)g(z)\,m(dz)\,m(dy)ds]\\
       &&=C^{1}_n\int_{0}^{t}\int_{\R^d}p(s,x,w)m(dw)\int_{\R^d}\frac{2C(w,y)(\ln(1+F(w,y))-F(w,y))1_{\{|w-y|>1/l\}}}{|w-y|^{d+\alpha}}\\
       &&\textrm{ }\textrm{ }\cdot \int_{\R^d}q_{n-1}(t-s,y,z)g(z)m(dz)\, m(dy)ds +\sum_{i=2}^nC^{i}_n\int_{0}^{t}\int_{\R^d}p(s,x,w)m(dw)\\
       &&\textrm{ }\textrm{ }\cdot \int_{\R^d}\frac{2C(w,y)(\ln(1+F(w,y)))^{i}1_{\{|w-y|>1/l\}}}{|w-y|^{d+\alpha}}\cdot\int_{\R^d}q_{n-i}(t-s,y,z)g(z)\, m(dz)\,m(dy)ds.\\
 \end{eqnarray*}

Let $\bf{K_{d,\alpha}}$ be the Kato class. Recall the definition
of Kato class. We say that a function $V$ on $\R^d$ belongs to the
Kato class $\bf{K_{d,\alpha}}$ if $\lim_{t\downarrow 0}\sup_{x\in
\R^d}\int_{0}^{t}\int_{\R^d}p(t,x,y)|V(y)|dyds=0.$ We say that a
signed measure $\mu$ on $\R^d$ belongs to the Kato class
$\bf{K_{d,\alpha}}$ if $\lim_{t\downarrow 0}\sup_{x\in
\R^d}\int_{0}^{t}\int_{\R^d}p(t,x,y)|\mu|(dy)ds=0.$

Suppose $F$ is a function on $\R^d \times \R^d$.
\begin{defn} We
say $F$ belongs to $\bf{J_{d,\alpha}}$ if $F$ is bounded,
vanishing on the diagonal, and the function
$$
               x\mapsto \int_{\R^d}\frac{|F(x,y)|}{|x-y|^{d+\alpha}}\,dy
$$
belongs to $\bf{K_{d,\alpha}}$.
\end{defn}

We assume $\inf_{x,y\in E}F(x,y)>-1$, $F(x,y)$ vanishes on the
diagonal and satisfies
\begin{eqnarray*}
           && w\mapsto \int_{\R^d}\frac{F^2(w,y)}{|w-y|^{d+\alpha}}\,dy\\
           \textrm{and}\\
           && y\mapsto \int_{\R^d}\frac{F^2(w,y)}{|w-y|^{d+\alpha}}\,dw\\
\end{eqnarray*}
both belong to $\bf{K_{d,\alpha}}$, i.e. $F^2(x,y)$ and $F^2(y,x)$
both belongs to $\bf{J_{d,\alpha}}$.

It is clear that there exist constants $C_{1}$ and $C_{2}$  such
that
$$|\ln(1+F(x,y))-F(x,y)|\leq
C_{1}F^2(x,y),\,\,\forall x,y\in E,$$ and
$$|\ln(1+F(x,y))|\leq C_{2} |F(x,y)|,\,\,\forall x,y \in E.$$

For the constants $C_{1}$ and $C_{2}$ given above, we assume that
there exist constants $\overline{C}$, $L$ and $\overline{M}$ such
that
\begin{eqnarray}
\max(|2C(x,y)C_1|,|2C(x,y)C_2|,C_2)&\leq& \overline{C}\\
\max(|2F(x,y)|, |2\overline{C}F(x,y)|)&\leq& L,
\end{eqnarray}
and $0 < M(y)\leq \overline{M}$ where $m(dy)=M(y)dy.$

Define $\overline{F}(x,y)=|F(x,y)|+|F(y,x)|$, which is symmetric
and satisfies $|\overline{F}(w,y)| \le L$ . Define
$\overline{p}(t,x,y)=p(t,x,y)+p(t,y,x)$. Then
$\overline{p}(t,x,y)$ is symmetric and satisfies
$$2M_{1}t^{-\frac{d}{\alpha}}\left(1 \wedge \frac{t^{\frac{1}{\alpha}}}{|x-y|}\right)^{d+\alpha}\leq \overline{p}(t,x,y)\leq 2 M_{2}t^{-\frac{d}{\alpha}}\left(1 \wedge \frac{t^{\frac{1}{\alpha}}}{|x-y|}\right)^{d+\alpha},\,\,\forall (t,x,y)\in (0,\infty)\times\R^d\times\R^d. $$

Denote $
(\int_{\R^d}\frac{\overline{F}^2(w,y)}{|w-y|^{d+\alpha}}\,dy)dw$
by $\mu(dw)$ and let
$C_t=\sup_{x\in\R^d}\int_{0}^{t}\int_{\R^d}\overline{p}(s,x,w)\,\mu(dw)ds$.
Then $C_t \downarrow 0$ as $t \downarrow 0$. It is clear that
there exist two positive constants $D_{1}$ and $D_{2}$ such that
$D_{1}\le \int_{\R^d}\overline{p}(t,x,y)\,m(dy)\le D_{2}$, as
$\overline{p}(t,x,y)$ is comparable to $p(t,x,y)$. Let
$\overline{q}_{0}(t,x,z)=\overline{p}(s,x,z)$ and define
$\overline{q}_{n}(t,x,z)$ by
\begin{eqnarray*}
               \overline{q}_n(t,x,z)&=&C^{1}_n\int_{0}^{t}\int_{\R^d}\overline{p}(s,x,w)m(dw)\int_{\R^d}\frac{\overline{C}\overline{F}^{2}(w,y)}{|w-y|^{d+\alpha}}\overline{q}_{n-1}(t-s,y,z)\,m(dy)ds\\
               &&+\sum^n_{i=2}C^{i}_n\int_{0}^{t}\int_{\R^d}\overline{p}(s,x,w)m(dw)\int_{\R^d}\frac{\overline{C}^i\overline{F}^{i}(w,y)}{|w-y|^{d+\alpha}}\overline{q}_{n-i}(t-s,y,z)\,m(dy)ds.\\
\end{eqnarray*}
By the choice of $\overline{C}$ in (3.1), it is clear that
$|q_{n}(t,x,z)| \leq \overline{q}_{n}(t,x,z)$. We can see that
$q_n(t,x,z)$ may not be symmetric even when $F(x,z)$ and
$p(t,x,z)$ are symmetric. Thus we can not use an argument similar
as that used in lemma 3.3 in \cite{RW}.

Next are the main results.

For the positive constants $L,\,\overline{C},\,\textrm{ and
}\overline{M}$ given above, and $K<1,\,M_{2}$ given in Lemma 3.1,
Lemma 3.2 and (1.1), there exists positive constants $\tilde{C}$
and $\tilde{C}_{1}$ such that
\begin{eqnarray*}
      \overline{C}\overline{M}^2D_2&\leq& \tilde{C}K,\\
      L^{m-2}\overline{C}^2\overline{M}^2 &\leq& \frac{1}{2}\tilde{C}m!K^m,\,\, \forall m>0,\\
      \overline{C}\frac{L^2D^2_{2}}{(1-\frac{1}{10}-2^{-\frac{1}{2}})^{d+\alpha}}&\leq& \frac{1}{4}\tilde{C}_{1}K,\,\, \forall m>0\\
      \overline{C}\frac{L^mD^2_{2}}{(1-\frac{1}{10}-2^{-\frac{1}{2}})^{d+\alpha}}&\leq&\frac{1}{4}\tilde{C}_1m!K^m, \,\, \forall m>0.\\
\end{eqnarray*}

We claim
\begin{thm}
There exist constants $t_{2}>0,\,\tilde{C}>0\textrm { and } \tilde
{C}_1>0 $ such that when $0<t \le t_{2}$,
\begin{eqnarray}
     &1.&\textrm{ } \int_{\R^d}\overline{q}_{n}(t,x,z)m(dz)\le
      \tilde{C}C_{t}n!K^n, \,\, \forall n \ge 1,\\
      &{  }&\textrm{ }\int_{\R^d}\overline{q}_{n}(t,x,z)m(dx)\le
      \tilde{C}C_{t}n!K^n, \,\, \forall n \ge 1.\\
      &2.&\textrm{ }  \overline{q}_{n}(t,x,z) \le
\tilde{C}_{1}n!K^{n}t^{-\frac{d}{\alpha}}\left(1 \wedge
\frac{t^{\frac{1}{\alpha}}}{|x-z|}\right)^{d+\alpha},\,\,\forall
n\ge 0
\end{eqnarray}
\end{thm}
\pf We show statement $1$ holds for $n=1$
\begin{eqnarray*}
      &&\int_{\R^d}\overline{q}_{1}(t,x,z)\,m(dz)\\
      &&\leq \int_{\R^d}\int_{0}^{t}\int_{\R^d}\overline{p}(s,x,w)m(dw)\int_{\R^d}\frac{\overline{C}\overline{F}^2(w,y)}{|w-y|^{d+\alpha}}\overline{p}(t-s,y,z)\,m(dy)ds\,m(dz)\\
      &&=\int_{0}^{t}\int_{\R^d}\overline{p}(s,x,w)m(dw)\int_{\R^d}\frac{\overline{C}\overline{F}^2(w,y)}{|w-y|^{d+\alpha}}\int_{\R^d}\overline{p}(t-s,y,z)\,m(dz)\,m(dy)ds\\
      &&=\int_{0}^{t}\int_{\R^d}\overline{p}(s,x,w)m(dw)\int_{\R^d}\frac{\overline{C}\overline{F}^2(w,y)}{|w-y|^{d+\alpha}}\,m(dy)dsD_2\\
      &&\leq \overline{C}\overline{M}^2C_{t}D_2 \leq \tilde{C}C_{t}K\\
\end{eqnarray*}
and
\begin{eqnarray*}
      &&\int_{\R^d}\overline{q}_{1}(t,x,z)\,m(dx)\\
      &&\leq \int_{\R^d}\int_{0}^{t}\int_{\R^d}\overline{p}(s,x,w)m(dw)\int_{\R^d}\frac{\overline{C}\overline{F}^2(w,y)}{|w-y|^{d+\alpha}}\overline{p}(t-s,y,z)\,m(dy)ds\,m(dx)\\
\end{eqnarray*}
 \begin{eqnarray*}
      &&=\int_{0}^{t}\int_{\R^d}\int_{\R^d}\overline{p}(s,x,w)\,m(dx)\,m(dw)\int_{\R^d}\frac{\overline{C}\overline{F}^2(w,y)}{|w-y|^{d+\alpha}}\overline{p}(t-s,y,z)\,m(dy)ds\\
      &&=\int_{0}^{t}\int_{\R^d}\int_{\R^d}\overline{p}(s,w,x)\,m(dx)\,m(dw)\int_{\R^d}\frac{\overline{C}\overline{F}^2(w,y)}{|w-y|^{d+\alpha}}\overline{p}(t-s,y,z)\,m(dy)ds\\
      &&\textrm{ }\textrm{ } \,\,(\textrm{ by symmetry of } \overline{p}(s,x,w) )\\
      &&=\int_{0}^{t}\int_{\R^d}\int_{\R^d}\frac{\overline{C}\overline{F}^2(w,y)}{|w-y|^{d+\alpha}}\,m(dw)p(t-s,y,z)\,m(dy)ds\\
      &&\leq \overline{C}\overline{M}^2C_{t}D_2 \leq \tilde{C}C_{t}K.\\
\end{eqnarray*}
Thus statement $1$ holds for $n=1$.

Since $\overline{q}_{0}(t,x,z)=\overline{p}(t,x,z)$, statement $2$
holds for $n=0$. We show it also holds for $n=1$. We write
$\overline{q}_1(t,x,z)$ in the following way,
\begin{eqnarray*}
      \overline{q}_1(t,x,z)&=&\int_{0}^{\frac{t}{2}}\int_{\R^d}\overline{p}(s,x,w)m(dw)\int_{\R^d}\frac{\overline{C}\overline{F}^2(w,y)}{|w-y|^{d+\alpha}}\overline{p}(t-s,y,z)\,m(dy)ds\\
      &&+\int_{\frac{t}{2}}^{t}\int_{\R^d}\overline{p}(s,x,w)m(dw)\int_{\R^d}\frac{\overline{C}\overline{F}^2(w,y)}{|w-y|^{d+\alpha}}\overline{p}(t-s,y,z)\,m(dy)ds.\\
\end{eqnarray*}

First we look at the first term in the expression of
$\overline{q}_1(t,x,z)$.

There are two cases:

Case 1. When $|x-z|\le t^{\frac{1}{\alpha}}$,
\begin{eqnarray*}
       &&\int_{0}^{\frac{t}{2}}\int_{\R^d}\overline{p}(s,x,w)m(dw)\int_{\R^d}\frac{\overline{C}\overline{F}^2(w,y)}{|w-y|^{d+\alpha}}\overline{p}(t-s,y,z)\,m(dy)ds\\
       &&\leq \int_{0}^{\frac{t}{2}}\int_{\R^d}\overline{p}(s,x,w)m(dw)\int_{\R^d}\frac{\overline{F}^2(w,y)}{|w-y|^{d+\alpha}}\,dyds\overline{C}\overline{M}\tilde{C}_{1}(\frac{t}{2})^{-\frac{d}{\alpha}}\\
       &&\leq C_t\overline{C}{\overline{M}}^{2}\tilde{C}_{1}(\frac{t}{2})^{-\frac{d}{\alpha}}.\\
\end{eqnarray*}
Since $C_t\downarrow 0$ as $t\downarrow 0$, there exists
$t_{01}>0$ such that when $0<t\leq t_{01}$,
$$C_t\overline{C}{\overline{M}}^{2}(\frac{1}{2})^{-\frac{d}{\alpha}}\leq
\frac{1}{2}K.$$ Thus
$$\int_{0}^{\frac{t}{2}}\int_{\R^d}\overline{p}(s,x,w)m(dw)\int_{\R^d}\frac{\overline{C}\overline{F}^2(w,y)}{|w-y|^{d+\alpha}}\overline{p}(t-s,y,z)\,m(dy)ds\leq
\frac{1}{2}\tilde{C}_{1}Kt^{-\frac{d}{\alpha}}.
$$

Case 2. When $|x-z|\ge t^{\frac{1}{\alpha}}$. Let
$B_1=\{y\in\R^d|\,|y-z| \ge \frac{1}{10}|x-z|\},
B_2=\{w\in\R^d|\,|w-x| \ge  2^{-\frac{1}{2}}|x-z|\}$ and
$B_3=\{(w,y)\in \R^d \times \R^d|\,|y-z| <  \frac{1}{10}|x-z|,\,
|w-x| <  2^{-\frac{1}{2}}|x-z|\} $. On $B_3$, we have $|w-y| \ge
(1-\frac{1}{10}- 2^{-\frac{1}{2}} )|x-z|$.
\begin{eqnarray*}
      &&\int_{0}^{\frac{t}{2}}\int_{\R^d}\overline{p}(s,x,w)m(dw)\int_{\R^d}\frac{\overline{C}\overline{F}^2(w,y)}{|w-y|^{d+\alpha}}\overline{p}(t-s,y,z)\,m(dy)ds\\
      &&\leq \int_{0}^{\frac{t}{2}}\int_{\R^d}\overline{p}(s,x,w)m(dw)\int_{\R^d}\frac{\overline{C}\overline{F}^2(w,y)}{|w-y|^{d+\alpha}}\overline{p}(t-s,y,z)1_{B_1}(y)\,m(dy)ds\\
      &&\textrm{  }\textrm{  }+\int_{0}^{\frac{t}{2}}\int_{\R^d}\overline{p}(s,x,w)m(dw)\int_{\R^d}\frac{\overline{C}\overline{F}^2(w,y)}{|w-y|^{d+\alpha}}\overline{p}(t-s,y,z)1_{B_2}(w)\,m(dy)ds\\
      &&\textrm{  }\textrm{  }+\int_{0}^{\frac{t}{2}}\int_{\R^d}\overline{p}(s,x,w)m(dw)\int_{\R^d}\frac{\overline{C}\overline{F}^2(w,y)}{|w-y|^{d+\alpha}}\overline{p}(t-s,y,z)1_{B_3}(w,y)\,m(dy)ds\\
      &&\leq \int_{0}^{\frac{t}{2}}\int_{\R^d}\overline{p}(s,x,w)m(dw)\int_{\R^d}\frac{\overline{F}^2(w,y)}{|w-y|^{d+\alpha}}\overline{M}\overline{C}\tilde{C}_{1}10^{d+\alpha}\frac{(t-s)}{|x-z|^{d+\alpha}}\,dyds\\
      &&\textrm{ }\textrm{ } +\int_{0}^{\frac{t}{2}}\int_{\R^d}dw\tilde{C}_{1}{\overline{M}}^22^{\frac{1}{2}(d+\alpha)}\frac{s}{{|x-z|}^{d+\alpha}}\int_{\R^d}\frac{\overline{F}^2(w,y)}{|w-y|^{d+\alpha}}\overline{C}\overline{p}(t-s,y,z)\,dyds\\
      &&\textrm{ }\textrm{ } +\overline{C}\frac{L^2}{(1-\frac{1}{10}- 2^{-\frac{1}{2}})^{d+\alpha}}\int_{0}^{\frac{t}{2}}\int_{\R^d}\overline{p}(s,x,w)m(dw)\int_{\R^d}\frac{1}{|x-z|^{d+\alpha}}\overline{p}(t-s,y,z)\,m(dy)ds\\
      &&\leq \left(\tilde{C}_{1}(\overline{C}\overline{M}^210^{d+\alpha})C_t
       +\tilde{C}_{1}(\overline{C}\overline{M}^22^{\frac{1}{2}(d+\alpha)})C_t\right)\frac{t}{{|x-z|}^{d+\alpha}}
       +\overline{C}\frac{L^2D^2_{2}}{(1-\frac{1}{10}-2^{-\frac{1}{2}})^{d+\alpha}}\frac{t}{{|x-z|}^{d+\alpha}}\\
      &&\textrm{   }\textrm{   } \,\,(\textrm{by symmetry, }\, \overline{p}(t-s,y,z)=\overline{p}(t-s,z,y)  \textrm{ }).\\
\end{eqnarray*}
There exists $ t_{02}>0$ such that when $0<t\leq t_{02}$, $C_t$ is
small enough and the sum of the first two terms $\leq
\frac{1}{4}\tilde{C}_{1}K\frac{t}{{|x-z|}^{d+\alpha}}$. For the
third term, we know that
$$\overline{C}\frac{L^2D^2_{2}}{(1-\frac{1}{10}-2^{-\frac{1}{2}})^{d+\alpha}}\leq
\frac{1}{4}\tilde{C}_{1}K.$$ Thus
      $$\int_{0}^{\frac{t}{2}}\int_{\R^d}\overline{p}(s,x,w)m(dw)\int_{\R^d}\frac{\overline{C}\overline{F}^2(w,y)}{|w-y|^{d+\alpha}}\overline{p}(t-s,y,z)\,m(dy)ds\leq
      \frac{1}{2}\tilde{C}K\frac{t}{{|x-z|}^{d+\alpha}}.$$
Combining case 1 and case 2, when $0<t\leq \min(t_{01},t_{02})$,
       $$\int_{0}^{\frac{t}{2}}\int_{\R^d}\overline{p}(s,x,w)m(dw)\int_{\R^d}\frac{\overline{C}\overline{F}^2(w,y)}{|w-y|^{d+\alpha}}\overline{p}(t-s,y,z)\,m(dy)ds\leq \frac{1}{2}\tilde{C}_{1}Kt^{-\frac{d}{\alpha}}\left(1
       \wedge\frac{t^{\frac{1}{\alpha}}}{|x-z|}\right)^{d+\alpha}.$$

For the second term in the expression of
$\overline{q}_{1}(t,x,z)$:
$$\int_{\frac{t}{2}}^{t}\int_{\R^d}\overline{p}(s,x,w)m(dw)\int_{\R^d}\frac{\overline{C}\overline{F}^2(w,y)}{|w-y|^{d+\alpha}}\overline{p}(t-s,y,z)\,m(dy)ds,$$
let $t-s=\tilde{s}$. The second term becomes
$$\int_{0}^{\frac{t}{2}}\int_{\R^d}\overline{p}(t-\tilde{s},x,w)m(dw)\int_{\R^d}\frac{\overline{C}\overline{F}^2(w,y)}{|w-y|^{d+\alpha}}\overline{p}(\tilde{s},y,z)\,m(dy)d\tilde{s}.$$

There are two cases:

Case a. When $|x-z| \le t^{\frac{1}{\alpha}}$,
\begin{eqnarray*}
      &&\int_{0}^{\frac{t}{2}}\int_{\R^d}\overline{p}(t-\tilde{s},x,w)m(dw)\int_{\R^d}\frac{\overline{C}\overline{F}^2(w,y)}{|w-y|^{d+\alpha}}\overline{p}(\tilde{s},y,z)\,m(dy)d\tilde{s}\\
      &&\leq \int_{0}^{\frac{t}{2}}\int_{\R^d}\tilde{C}_{1}(\frac{t}{2})^{-\frac{d}{\alpha}}m(dw)\int_{\R^d}\frac{\overline{F}^2(w,y)}{|w-y|^{d+\alpha}}\overline{C}\overline{p}(\tilde{s},z,y)\,m(dy)d\tilde{s}\\
      &&\leq \tilde{C}_{1}\overline{C}\overline{M}^2(\frac{t}{2})^{-\frac{d}{\alpha}}C_{t}.\\
\end{eqnarray*}
Since $C_t\downarrow 0$ as $t\downarrow 0$, there exists $
t_{03}>0$ such that when $0<t\leq t_{03},$
$$\overline{C}\overline{M}^2(\frac{1}{2})^{-\frac{d}{\alpha}}C_t\leq
\frac{1}{2}K.$$
 Thus
 $$\int_{0}^{\frac{t}{2}}\int_{\R^d}\overline{p}(t-\tilde{s},x,w)m(dw)\int_{\R^d}\frac{\overline{C}\overline{F}^2(w,y)}{|w-y|^{d+\alpha}}\overline{p}(\tilde{s},y,z)\,m(dy)d\tilde{s}\leq
 \frac{1}{2}\tilde{C}_{1}Kt^{-\frac{d}{\alpha}}.$$

Case b. When $|x-z|\ge t^{\frac{1}{\alpha}}$. Let
$B_1=\{y\in\R^d|\,|y-z| \ge \frac{1}{10}|x-z|\},
B_2=\{w\in\R^d|\,|w-x| \ge  2^{-\frac{1}{2}}|x-z|\}$ and
$B_3=\{(w,y)\in \R^d \times \R^d|\,|y-z| <  \frac{1}{10}|x-z|,\,
|w-x| <  2^{-\frac{1}{2}}|x-z|\} $. On $B_3$, we have $|w-y| \ge
(1-\frac{1}{10}- 2^{-\frac{1}{2}} )|x-z|$.
\begin{eqnarray*}
      &&\int_{0}^{\frac{t}{2}}\int_{\R^d}\overline{p}(t-\tilde{s},x,w)m(dw)\int_{\R^d}\frac{\overline{C}\overline{F}^2(w,y)}{|w-y|^{d+\alpha}}\overline{p}(\tilde{s},y,z)\,m(dy)d\tilde{s}\\
      &&\leq \int_{0}^{\frac{t}{2}}\int_{\R^d}\overline{p}(t-\tilde{s},x,w)m(dw)\int_{\R^d}\frac{\overline{C}\overline{F}^2(w,y)}{|w-y|^{d+\alpha}}\overline{p}(\tilde{s},y,z)1_{B_1}(y)\,m(dy)d\tilde{s}\\
      &&\textrm{ }\textrm{ }+ \int_{0}^{\frac{t}{2}}\int_{\R^d}\overline{p}(t-\tilde{s},x,w)m(dw)\int_{\R^d}\frac{\overline{C}\overline{F}^2(w,y)}{|w-y|^{d+\alpha}}\overline{p}(\tilde{s},y,z)1_{B_2}(w)\,m(dy)d\tilde{s}\\
      &&\textrm{ }\textrm{ }+ \int_{0}^{\frac{t}{2}}\int_{\R^d}\overline{p}(t-\tilde{s},x,w)m(dw)\int_{\R^d}\frac{\overline{C}\overline{F}^2(w,y)}{|w-y|^{d+\alpha}}\overline{p}(\tilde{s},y,z)1_{B_3}(w,y)\,m(dy)d\tilde{s}\\
      &&\leq \int_{0}^{\frac{t}{2}}\int_{\R^d}\overline{p}(t-\tilde{s},x,w)m(dw)\int_{\R^d}\frac{\overline{F}^2(w,y)}{|w-y|^{d+\alpha}}\overline{C}\tilde{C}_{1}10^{d+\alpha}\frac{\tilde{s}}{|x-z|^{d+\alpha}}\,m(dy)d\tilde{s}\\
      &&\textrm{ }\textrm{ }+ \int_{0}^{\frac{t}{2}}\int_{\R^d}\,dw\tilde{C}_{1}\overline{M}^22^{\frac{1}{2}(d+\alpha)}\frac{t-\tilde{s}}{{|x-z|}^{d+\alpha}}\int_{\R^d}\frac{\overline{F}^2(w,y)}{|w-y|^{d+\alpha}}\overline{C}\overline{p}(\tilde{s},z,y)\,dyd\tilde{s}\\
\end{eqnarray*}
\begin{eqnarray*}
      &&\textrm{ }\textrm{ } +\overline{C}\frac{L^2}{(1-\frac{1}{10}- 2^{-\frac{1}{2}})^{d+\alpha}}\int_{0}^{\frac{t}{2}}\int_{\R^d}\overline{p}(t-\tilde{s},x,w)m(dw)\int_{\R^d}\frac{1}{|x-z|^{d+\alpha}}\overline{p}(\tilde{s},y,z)\,m(dy)d\tilde{s}\\
      &&\leq (\tilde{C}_{1}(\overline{C}\overline{M}^210^{d+\alpha})C_t+\tilde{C}_{1}(\overline{C}\overline{M}^22^{\frac{1}{2}(d+\alpha)})C_t\frac{t}{{|x-z|}^{d+\alpha}})+\overline{C}\frac{L^2D^2_2}{(1-\frac{1}{10}-2^{-\frac{1}{2}})^{d+\alpha}}\frac{t}{{|x-z|}^{d+\alpha}}.\\
\end{eqnarray*}
 There exists $ t_{04}>0$ such that when $0<t\leq t_{04}$, $C_t$ is small enough and
the sum of the first two terms $\leq
\frac{1}{4}\tilde{C}_{1}K\frac{t}{{|x-z|}^{d+\alpha}}$. For the
third term, we know that
$$\overline{C}\frac{L^2D^2_2}{(1-\frac{1}{10}-2^{-\frac{1}{2}})^{d+\alpha}}\leq
\frac{1}{4}\tilde{C}_{1}K.$$ Thus
$$\int_{0}^{\frac{t}{2}}\int_{\R^d}\overline{p}(t-\tilde{s},x,w)m(dw)\int_{\R^d}\frac{\overline{C}\overline{F}^2(w,y)}{|w-y|^{d+\alpha}}\overline{p}(\tilde{s},y,z)\,m(dy)d\tilde{s}\leq \frac{1}{2}\tilde{C}_{1}K\frac{t}{{|x-z|}^{d+\alpha}}.$$
Combining case a and case b, when $0<t\leq \min(t_{03},t_{04})$,
$$\int_{\frac{t}{2}}^{t}\int_{\R^d}\overline{p}(s,x,w)m(dw)\int_{\R^d}\frac{\overline{C}\overline{F}^2(w,y)}{|w-y|^{d+\alpha}}\overline{p}(t-s,y,z)\,m(dy)ds\leq \frac{1}{2}\tilde{C}_{1}Kt^{-\frac{d}{\alpha}}\left(1
\wedge\frac{t^{\frac{1}{\alpha}}}{|x-z|}\right)^{d+\alpha}.$$
Therefore, when $0<t\leq t_{0}=\min(t_{01},t_{02},t_{03},t_{04})$,
$$\overline{q}_{1}(t,x,z)\leq \tilde{C}_{1}Kt^{-\frac{d}{\alpha}}\left(1\wedge\frac{t^{\frac{1}{\alpha}}}{|x-z|}\right)^{d+\alpha}.$$
i.e. statement $2$ holds for $n=1$.

Next we suppose statement $1$ and  statement $2$ are true for
$n\leq m-1$, we show that they hold for $n=m$.

First we look at statement $1$.

For the first part of statement $1$,
\begin{eqnarray*}
      &&\int_{\R^d}\overline{q}_m(t,x,z)\,m(dz)\\
      &&\leq\int_{\R^d} C^{1}_m\int_{0}^{t}\int_{\R^d}\overline{p}(s,x,w)m(dw)\int_{\R^d}\frac{\overline{C}\overline{F}^2(w,y)}{|w-y|^{d+\alpha}}\overline{q}_{m-1}(t-s,y,z)\,m(dy)ds\,m(dz)\\
      &&\textrm{ }\textrm{ } +\int_{\R^d}\sum_{i=2}^mC^{i}_m\int_{0}^{t}\int_{\R^d}\overline{p}(s,x,w)m(dw)\int_{\R^d}\frac{\overline{C}^2\overline{F}^{2}(w,y)}{|w-y|^{d+\alpha}}L^{i-2}\overline{q}_{m-i}(t-s,y,z)\,m(dy)ds\,m(dz) \\
      &&\leq C^{1}_m\int_{0}^{t}\int_{\R^d}\overline{p}(s,x,w)m(dw)\int_{\R^d}\frac{\overline{C}\overline{F}^2(w,y)}{|w-y|^{d+\alpha}}(\int_{\R^d}\overline{q}_{m-1}(t-s,y,z)\,m(dz))\,m(dy)ds\\
\end{eqnarray*}
\begin{eqnarray*}
      &&\textrm{ }\textrm{ } +\sum_{i=2}^{m-1}C^{i}_{m}\int_{0}^{t}\int_{\R^d}\overline{p}(s,x,w)m(dw)\int_{\R^d}\frac{\overline{C}^2\overline{F}^{2}(w,y)}{|w-y|^{d+\alpha}}L^{i-2}(\int_{\R^d}\overline{q}_{m-i}(t-s,y,z)\,m(dz))\,m(dy)ds \\
      &&\textrm{ }\textrm{ } +\int_{0}^{t}\int_{\R^d}\overline{p}(s,x,w)m(dw)\int_{\R^d}\frac{\overline{C}^2\overline{F}^{2}(w,y)}{|w-y|^{d+\alpha}}L^{m-2}(\int_{\R^d}\overline{p}(t-s,y,z)\,m(dz))\,m(dy)ds \\
      &&\leq C^{1}_m\tilde{C}C_{t}(m-1)!K^{m-1}\int_{0}^{t}\int_{\R^d}\overline{p}(s,x,w)m(dw)\int_{\R^d}\frac{\overline{C}\overline{F}^2(w,y)}{|w-y|^{d+\alpha}}\,m(dy)ds\\
      &&\textrm{ }\textrm{ } +\sum_{i=2}^{m-1}C^{i}_{m}\tilde{C}C_{t}(m-i)!L^{i-2}K^{m-i}\int_{0}^{t}\int_{\R^d}\overline{p}(s,x,w)m(dw)\int_{\R^d}\frac{\overline{C}^2\overline{F}^{2}(w,y)}{|w-y|^{d+\alpha}}\,m(dy)ds \\
      &&\textrm{ }\textrm{ } \,\,(\textrm{ since statement 1 holds for } n\leq m-1 )\\
      &&\textrm{ }\textrm{ } +\int_{0}^{t}\int_{\R^d}\overline{p}(s,x,w)m(dw)\int_{\R^d}\frac{\overline{C}^2\overline{F}^{2}(w,y)}{|w-y|^{d+\alpha}}L^{m-2}\,m(dy)dsD_2 \\
      &&\leq \tilde{C}C_{t}(m-1)!K^{m-1}\overline{C}\overline{M}^2C_{t}+\sum_{i=2}^{m-1}C^{i}_{m}\tilde{C}C_{t}(m-i)!L^{i-2}K^{m-i}\overline{C}^2\overline{M}^2C_{t}
       +L^{m-2}\overline{C}^2\overline{M}^2D_2C_{t}.\\
\end{eqnarray*}
For the sum of the first and the second term, there exists $
t_{11}>0 $ with $t_{11}\leq t_{0}$ such that when $0<t< t_{11}$,
$C_t$ is small enough and
      $$\tilde{C}C_{t}\left(C^{1}_m(m-1)!K^{m-1}+\sum_{i=2}^{m-1}C^{i}_{m}(m-i)!L^{i-2}K^{m-i}\overline{C}\right)\overline{C}\overline{M}^2\leq \frac{1}{2}\tilde{C}m!K^m.$$
For the third term, we know that
      $$L^{m-2}\overline{C}^2\overline{M}^2D_2\leq \frac{1}{2}\tilde{C}m!K^m,\,\, \forall m>0.$$
Thus
      $$\int_{\R^d}\overline{q}_m(t,x,z)\,m(dz)\leq \tilde{C}C_{t}m!K^m.$$
i.e. the first part of the statement $1$ holds for $n=m$.

To show the second part of the statement $1$ for $n=m$, we know
that
\begin{eqnarray*}
&&\overline{q}_m(t,x,z)\\
 &&\leq C^{1}_m\int_{0}^{t}\int_{\R^d}\overline{p}(s_{1},x,w_{1})m(dw_{1})\int_{\R^d}\frac{\overline{C}\overline{F}^2(w_{1},y_{1})}{|w_{1}-y_{1}|^{d+\alpha}}\overline{q}_{m-1}(t-s_{1},y_{1},z)\,m(dy_{1})ds_{1}\\
      &&\textrm{ }\textrm{ }+\sum_{i=2}^mC^{i}_m\int_{0}^{t}\int_{\R^d}\overline{p}(s_{1},x,w_{1})m(dw_{1})\int_{\R^d}\frac{\overline{C}^2\overline{F}^{2}(w_{1},y_{1})}{|w_{1}-y_{1}|^{d+\alpha}}L^{i-2}\overline{q}_{m-i}(t-s_{1},y_{1},z)\,m(dy_{1})ds_{1}. \\
\end{eqnarray*}
Substituting the expressions of
$\overline{q}_{m-1}(t-s_{1},y_{1},z)$ and
$\overline{q}_{m-i}(t-s_{1},y_{1},z)$ into the above, we have
\begin{eqnarray*}
      &&\overline{q}_m(t,x,z)\\
      &&\leq C^{1}_m\int_{0}^{t}\int_{\R^d}\overline{p}(s_{1},x,w_{1})m(dw_{1})\int_{\R^d}\frac{\overline{C}\overline{F}^2(w_{1},y_{1})}{|w_{1}-y_{1}|^{d+\alpha}}\,m(dy_{1})ds_{1}(C^{1}_{m-1}\int_{0}^{t-s_1}\int_{\R^d}\overline{p}(s_{2},y_1,w_{2})\\
      &&\textrm{ }\textrm{ }\cdot m(dw_{2})\int_{\R^d}\frac{\overline{C}\overline{F}^2(w_{2},y_{2})}{|w_{2}-y_{2}|^{d+\alpha}}\overline{q}_{m-2}(t-s_{1}-s_2,y_{2},z)\,m(dy_{2})ds_{2})\\
      &&\textrm{ }\textrm{ }+C^{1}_m\int_{0}^{t}\int_{\R^d}\overline{p}(s_{1},x,w_{1})m(dw_{1})\int_{\R^d}\frac{\overline{C}\overline{F}^2(w_{1},y_{1})}{|w_{1}-y_{1}|^{d+\alpha}}\,m(dy_{1})ds_{1}(\sum_{i=2}^{m-1}C^{i}_{m-1}\int_{0}^{t-s_1}\\
      &&\textrm{ }\textrm{ }\cdot\int_{\R^d}\overline{p}(s_{2},y_1,w_{2})m(dw_{2})\int_{\R^d}\frac{\overline{C}^2\overline{F}^2(w_{2},y_{2})}{|w_{2}-y_{2}|^{d+\alpha}}L^{i-2}\overline{q}_{m-1-i}(t-s_{1}-s_2,y_{2},z)\,m(dy_{2})ds_{2})\\
      &&\textrm{ }\textrm{ } +\sum_{i=2}^mC^{i}_m\int_{0}^{t}\int_{\R^d}\overline{p}(s_{1},x,w_{1})m(dw_{1})\int_{\R^d}\frac{\overline{C}^2\overline{F}^2(w_{1},y_{1})}{|w_{1}-y_{1}|^{d+\alpha}}L^{i-2}\,m(dy_{1})ds_{1}(C^{1}_{m-i}\int_{0}^{t-s_1}\\
      &&\textrm{ }\textrm{ } \cdot \int_{\R^d}\overline{p}(s_{2},y_1,w_{2})m(dw_{2})\int_{\R^d}\frac{\overline{C}\overline{F}^2(w_{2},y_{2})}{|w_{2}-y_{2}|^{d+\alpha}}\overline{q}_{m-i-1}(t-s_{1}-s_2,y_{2},z)\,m(dy_{2})ds_{2})\\
      &&\textrm{ }\textrm{ } +\sum_{i=2}^mC^{i}_m\int_{0}^{t}\int_{\R^d}\overline{p}(s_{1},x,w_{1})m(dw_{1})\int_{\R^d}\frac{\overline{C}^2\overline{F}^2(w_{1},y_{1})}{|w_{1}-y_{1}|^{d+\alpha}}L^{i-2}\,m(dy_{1})ds_{1}(\sum_{j=2}^{m-i}C^{j}_{m-i}\int_{0}^{t-s_1}\\
      &&\textrm{ }\textrm{ } \cdot \int_{\R^d}\overline{p}(s_{2},y_1,w_{2})m(dw_{2})\int_{\R^d}\frac{\overline{C}^2\overline{F}^2(w_{2},y_{2})}{|w_{2}-y_{2}|^{d+\alpha}}L^{j-2}\overline{q}_{m-i-j}(t-s_{1}-s_2,y_{2},z)\,m(dy_{2})ds_{2})\\
      &&\leq C^{1}_m\int_{0}^{t}\int_{\R^d}\overline{p}(s_{1},x,w_{1})m(dw_{1})\int_{\R^d}\frac{\overline{C}\overline{F}^2(w_{1},y_{1})}{|w_{1}-y_{1}|^{d+\alpha}}\,m(dy_{1})ds_{1}(C^{1}_{m-1}\int_{0}^{t-s_1}\int_{\R^d}\overline{p}(s_{2},y_1,w_{2})\\
      &&\textrm{ }\textrm{ }\cdot m(dw_{2})\int_{\R^d}\frac{\overline{C}\overline{F}^2(w_{2},y_{2})}{|w_{2}-y_{2}|^{d+\alpha}}\tilde{C}_{1}(m-2)!K^{m-2}\frac{1}{2M_1}\overline{p}(t-s_{1}-s_2,y_{2},z)\,m(dy_{2})ds_{2})\\
      &&\textrm{ }\textrm{ }+C^{1}_m\int_{0}^{t}\int_{\R^d}\overline{p}(s_{1},x,w_{1})m(dw_{1})\int_{\R^d}\frac{\overline{C}\overline{F}^2(w_{1},y_{1})}{|w_{1}-y_{1}|^{d+\alpha}}\,m(dy_{1})ds_{1}(\sum_{i=2}^{m-1}C^{i}_{m-1}\int_{0}^{t-s_1}\int_{\R^d}\\
      &&\textrm{ }\textrm{ }\cdot \overline{p}(s_{2},y_1,w_{2})m(dw_{2})\int_{\R^d}\frac{\overline{C}^2\overline{F}^2(w_{2},y_{2})}{|w_{2}-y_{2}|^{d+\alpha}}\tilde{C}_{1}(m-1-i)!L^{i-2}K^{m-1-i}\frac{1}{2M_1}\\
      &&\textrm{ }\textrm{ }\cdot \overline{p}(t-s_{1}-s_2,y_{2},z)\,m(dy_{2})ds_{2})\\
      &&\textrm{ }\textrm{ } +\sum_{i=2}^mC^{i}_m\int_{0}^{t}\int_{\R^d}\overline{p}(s_{1},x,w_{1})m(dw_{1})\int_{\R^d}\frac{\overline{C}^2\overline{F}^2(w_{1},y_{1})}{|w_{1}-y_{1}|^{d+\alpha}}L^{i-2}\,m(dy_{1})ds_{1}(C^{1}_{m-i}\int_{0}^{t-s_1}\int_{\R^d}\\
      &&\textrm{ }\textrm{ } \cdot \overline{p}(s_{2},y_1,w_{2})m(dw_{2})\int_{\R^d}\frac{\overline{C}\overline{F}^2(w_{2},y_{2})}{|w_{2}-y_{2}|^{d+\alpha}}\tilde{C}_{1}(m-i-1)!K^{m-i-1}\frac{1}{2M_1}\\
      &&\textrm{ }\textrm{ }\cdot \overline{p}(t-s_{1}-s_2,y_{2},z)\,m(dy_{2})ds_{2})\\
      &&\textrm{ }\textrm{ } +\sum_{i=2}^mC^{i}_m\int_{0}^{t}\int_{\R^d}\overline{p}(s_{1},x,w_{1})m(dw_{1})\int_{\R^d}\frac{\overline{C}^2\overline{F}^2(w_{1},y_{1})}{|w_{1}-y_{1}|^{d+\alpha}}L^{i-2}\,m(dy_{1})ds_{1}(\sum_{j=2}^{m-i}C^{j}_{m-i}\int_{0}^{t-s_1}\\
\end{eqnarray*}
\begin{eqnarray*}
      &&\textrm{ }\textrm{ } \cdot \int_{\R^d}\overline{p}(s_{2},y_1,w_{2}) m(dw_{2})\int_{\R^d}\frac{\overline{C}^2\overline{F}^2(w_{2},y_{2})}{|w_{2}-y_{2}|^{d+\alpha}}\tilde{C}_{1}(m-i-j)!L^{j-2}K^{m-i-j}\frac{1}{2M_1}\\
       &&\textrm{ }\textrm{ } \cdot\overline{p}(t-s_{1}-s_2,y_{2},z)\,m(dy_{2})ds_{2})\\
      &&\textrm{ }\textrm{ }\,\, (\textrm{ since the statement 2 holds true for } n\leq m-1 )\\
      &&= C^{1}_mC^{1}_{m-1}\overline{C}^2\tilde{C}_{1}(m-2)!K^{m-2}\frac{1}{2M_1}(\int_{0}^{t}\int_{\R^d}\overline{p}(s_{1},x,w_{1})m(dw_{1})\int_{\R^d}\frac{\overline{F}^2(w_{1},y_{1})}{|w_{1}-y_{1}|^{d+\alpha}}\\
      &&\textrm{ }\textrm{ }\cdot m(dy_{1})ds_{1}\int_{0}^{t-s_1}\int_{\R^d}\overline{p}(s_{2},y_1,w_{2})m(dw_{2})\int_{\R^d}\frac{\overline{F}^2(w_{2},y_{2})}{|w_{2}-y_{2}|^{d+\alpha}}\overline{p}(t-s_{1}-s_2,y_{2},z)\,m(dy_{2})ds_{2})\\
      &&\textrm{ }\textrm{ }+C^{1}_m\sum_{i=2}^{m-1}C^{i}_{m-1}\overline{C}^3\tilde{C}_{1}(m-1-i)!L^{i-2}K^{m-1-i}\frac{1}{2M_1}(\int_{0}^{t}\int_{\R^d}\overline{p}(s_{1},x,w_{1})m(dw_{1})\\
      &&\textrm{ }\textrm{ }\cdot \int_{\R^d}\frac{\overline{F}^2(w_{1},y_{1})}{|w_{1}-y_{1}|^{d+\alpha}}m(dy_{1})ds_{1}\int_{0}^{t-s_1}\int_{\R^d}\overline{p}(s_{2},y_1,w_{2})m(dw_{2})\int_{\R^d}\frac{\overline{F}^2(w_{2},y_{2})}{|w_{2}-y_{2}|^{d+\alpha}}\\
      &&\textrm{ }\textrm{ }\cdot\overline{p}(t-s_{1}-s_2,y_{2},z)\,m(dy_{2})ds_{2})\\
      &&\textrm{ }\textrm{ } +\sum_{i=2}^mC^{i}_mC^{1}_{m-i}\overline{C}^3\tilde{C}_{1}(m-i-1)!L^{i-2}K^{m-i-1}\frac{1}{2M_1}(\int_{0}^{t}\int_{\R^d}\overline{p}(s_{1},x,w_{1})m(dw_{1})\\
      &&\textrm{ }\textrm{ }\cdot \int_{\R^d}\frac{\overline{F}^2(w_{1},y_{1})}{|w_{1}-y_{1}|^{d+\alpha}}m(dy_{1})ds_{1}\int_{0}^{t-s_1}\int_{\R^d}\overline{p}(s_{2},y_1,w_{2})m(dw_{2})\int_{\R^d}\frac{\overline{F}^2(w_{2},y_{2})}{|w_{2}-y_{2}|^{d+\alpha}}\\
      &&\textrm{ }\textrm{ }\cdot\overline{p}(t-s_{1}-s_2,y_{2},z)\,m(dy_{2})ds_{2})\\
      &&\textrm{ }\textrm{ } +\sum_{i=2}^mC^{i}_m\sum_{j=2}^{m-i}C^{j}_{m-i}\overline{C}^4\tilde{C}_{1}(m-i-j)!L^{i+j-4}K^{m-i-j}\frac{1}{2M_1}(\int_{0}^{t}\int_{\R^d}\overline{p}(s_{1},x,w_{1})m(dw_{1})\\
      &&\textrm{ }\textrm{ }\cdot \int_{\R^d}\frac{\overline{F}^2(w_{1},y_{1})}{|w_{1}-y_{1}|^{d+\alpha}} m(dy_{1})ds_{1}\int_{0}^{t-s_1}\int_{\R^d}\overline{p}(s_{2},y_1,w_{2})m(dw_{2})\int_{\R^d}\frac{\overline{F}^2(w_{2},y_{2})}{|w_{2}-y_{2}|^{d+\alpha}}\\
      &&\textrm{ }\textrm{ }\cdot \overline{p}(t-s_{1}-s_2,y_{2},z)\,m(dy_{2})ds_{2}).\\
\end{eqnarray*}
It is clear that in the above expression there is the following
common factor :
\begin{eqnarray*}
      &&\int_{0}^{t}\int_{\R^d}\overline{p}(s_{1},x,w_{1})m(dw_{1})\int_{\R^d}\frac{\overline{F}^2(w_{1},y_{1})}{|w_{1}-y_{1}|^{d+\alpha}}\,m(dy_{1})ds_{1}\int_{0}^{t-s_1}\int_{\R^d}\overline{p}(s_{2},y_1,w_{2})m(dw_{2})\\
      &&\cdot\int_{\R^d}\frac{\overline{F}^2(w_{2},y_{2})}{|w_{2}-y_{2}|^{d+\alpha}}\overline{p}(t-s_{1}-s_2,y_{2},z)\,m(dy_{2})ds_{2}.\\
\end{eqnarray*}
We denote it by $G(t,x,z)$.

Next we show $ G(t,x,z)=G(t,z,x).$

 Changing variables, let
$t-s_1-s_2=\tilde{s}_1,\,s_2=\tilde{s}_2$. The absolute value of
the Jacobian of the this transformation is $1$. Let
$y_2=\tilde{w}_1,\, y_1=\tilde{w_2},\,w_2=\tilde{y}_1,
w_1=\tilde{y}_2$.\\
Thus
\begin{eqnarray*}
G(t,x,z)&=&\int_{0}^{t}\int_{0}^{t-\tilde{s}_1}\int_{\R^d}\overline{p}(t-\tilde{s}_1-\tilde{s}_2,x,\tilde{y}_2)m(d\tilde{y}_2)\int_{\R^d}\frac{\overline{F}^{2}(\tilde{y}_2,\tilde{w_2})}{|\tilde{y}_2-\tilde{w_2}|^{d+\alpha}}\,m(d\tilde{w_2})\\
      &&\cdot(\int_{\R^d}\overline{p}(\tilde{s}_{2},\tilde{w_2},\tilde{y}_1)m(d\tilde{y}_1)\int_{\R^d}\frac{\overline{F}^2(\tilde{y}_1,\tilde{w}_1)}{|\tilde{y}_1-\tilde{w}_1|^{d+\alpha}}\overline{p}(\tilde{s}_1,\tilde{w}_1,z)\,m(d\tilde{w}_1)d\tilde{s}_{2}\,d\tilde{s}_{1}.\\
\end{eqnarray*}
Rearranging the components of the integrand and using the fact
that $\overline{F}(x,y)$ and $\overline{p}(t,x,y)$ are symmetric
in $x$ and $y$, we can see that the  expression of $G(t,x,z)$ is
equal to $G(t,z,x)$.

Therefore
 $$\int_{\R^d}G(t,x,z)\,m(dx)
      =\int_{\R^d}G(t,z,x)\,m(dx),$$
and both
\begin{eqnarray*}
      &&=\int_{0}^{t}\int_{\R^d}\overline{p}(s_{1},z,w_{1})m(dw_{1})\int_{\R^d}\frac{\overline{F}^2(w_{1},y_{1})}{|w_{1}-y_{1}|^{d+\alpha}}\,m(dy_{1})ds_{1}(\int_{0}^{t-s_1}\int_{\R^d}\overline{p}(s_{2},y_1,w_{2})m(dw_{2})\\
      &&\textrm{ }\textrm{ }\cdot(\int_{0}^{t-s_1}\int_{\R^d}\overline{p}(s_{2},y_1,w_{2})m(dw_{2})\int_{\R^d}\frac{\overline{F}^2(w_{2},y_{2})}{|w_{2}-y_{2}|^{d+\alpha}}\int_{\R^d}\overline{p}(t-s_{1}-s_2,y_{2},x)m(dx)\,m(dy_{2})ds_{2})\\
      &&=\int_{0}^{t}\int_{\R^d}\overline{p}(s_{1},z,w_{1})m(dw_{1})\int_{\R^d}\frac{\overline{F}^2(w_{1},y_{1})}{|w_{1}-y_{1}|^{d+\alpha}}D_2\,m(dy_{1})ds_{1}\\
      &&\textrm{ }\textrm{ }\cdot(\int_{0}^{t-s_1}\int_{\R^d}\overline{p}(s_{2},y_1,w_{2})m(dw_{2})\int_{\R^d}\frac{\overline{F}^2(w_{2},y_{2})}{|w_{2}-y_{2}|^{d+\alpha}}\,m(dy_{2})ds_{2})\\
      &&\leq \int_{0}^{t}\int_{\R^d}\overline{p}(s_{1},z,w_{1})m(dw_{1})\int_{\R^d}\frac{\overline{F}^2(w_{1},y_{1})}{|w_{1}-y_{1}|^{d+\alpha}}\,m(dy_{1})ds_{1}
      (\int_{0}^{t-s_1}\int_{\R^d}\overline{p}(s_{2},y_1,w_{2})\,\mu(dw_{2})ds_{2}D_2\overline{M}^2)\\
      &&\leq \int_{0}^{t}\int_{\R^d}\overline{p}(s_{1},z,w_{1})m(dw_{1})\int_{\R^d}\frac{\overline{F}^2(w_{1},y_{1})}{|w_{1}-y_{1}|^{d+\alpha}}\,m(dy_{1})ds_{1}C_{t}D_2\overline{M}^2\\
      &&\leq \int_{0}^{t}\int_{\R^d}\overline{p}(s_{1},z,w_{1})\,\mu{dw_{1}}\overline{M}^2C_{t}D_2\overline{M}^2
      \leq D_2\overline{M}^4C_{t}C_{t},\\
\end{eqnarray*}
which implies that
\begin{eqnarray*}
      &&\int_{\R^d}\overline{q}_m(t,x,z)\,m(dx)\\
      &&\leq C^{1}_mC^{1}_{m-1}\overline{C}^2\tilde{C}_{1}(m-2)!K^{m-2}\frac{1}{2M_1}D_2\overline{M}^4C_{t}C_{t}\\
      &&\textrm{ }\textrm{ }+C^{1}_m\sum_{i=2}^{m-1}C^{i}_{m-1}\overline{C}^3\tilde{C}_{1}(m-1-i)!L^{i-2}K^{m-1-i}\frac{1}{2M_1}D_2\overline{M}^4C_{t}C_{t}\\
      &&\textrm{ }\textrm{ }+\sum_{i=2}^mC^{i}_mC^{1}_{m-i}\overline{C}^3\tilde{C}_{1}(m-i-1)!L^{i-2}K^{m-i-1}\frac{1}{2M_1}D_2\overline{M}^4C_{t}C_{t}\\
\end{eqnarray*}
      $$\textrm{ }\textrm{
      }+\sum_{i=2}^mC^{i}_m\sum_{j=2}^{m-i}C^{j}_{m-i}\overline{C}^4\tilde{C}_{1}(m-i-j)!L^{i-2}L^{j-2}K^{m-i-j}\frac{1}{2M_1}D_2\overline{M}^4C_{t}C_{t}.$$

Since $C_t \downarrow 0$ as $t\downarrow 0$ and Lemma 3.1 and
Lemma 3.2 hold, it is easy to see that there exists $ t_{12}>0$
with $t_{12}\leq t_{0}$ such that when $0<t\leq t_{12}$, the first
three term in the above expression $\leq \frac{1}{2}
\tilde{C}C_{t}m!K^m$ and the fourth term $\leq \frac{1}{2}
\tilde{C}C_{t}m!K^m$. Thus
$$\int_{\R^d}\overline{q}_m(t,x,z)\,m(dx)\leq \tilde{C}C_{t}m!K^m,$$
i.e. the second part of statement $1$ holds for $n=m$. Thus
statement $1$ holds for $n=m$.

Next we show statement $2$ for $n=m$.
 We write $\overline{q}_m(t,x,z)$ in the following way
\begin{eqnarray*}
&&\overline{q}_m(t,x,z)\\
      &&=C^{1}_m\int_{0}^{\frac{t}{2}}\int_{\R^d}\overline{p}(s,x,w)m(dw)\int_{\R^d}\frac{\overline{C}\overline{F}^2(w,y)}{|w-y|^{d+\alpha}}\overline{q}_{m-1}(t-s,y,z)\,m(dy)ds\\
       &&\textrm{ }\textrm{ } +\sum_{i=2}^mC^{i}_m\int_{0}^{\frac{t}{2}}\int_{\R^d}\overline{p}(s,x,w)m(dw)\int_{\R^d}\frac{\overline{C}^i\overline{F}^{i}(w,y)}{|w-y|^{d+\alpha}}\overline{q}_{m-i}(t-s,y,z)\,m(dy)ds. \\
       &&\textrm{ }\textrm{ } +C^{1}_m\int_{\frac{t}{2}}^{t}\int_{\R^d}\overline{p}(s,x,w)m(dw)\int_{\R^d}\frac{\overline{C}\overline{F}^2(w,y)}{|w-y|^{d+\alpha}}\overline{q}_{m-1}(t-s,y,z)\,m(dy)ds\\
       &&\textrm{ }\textrm{ } +\sum_{i=2}^mC^{i}_m\int_{\frac{t}{2}}^{t}\int_{\R^d}\overline{p}(s,x,w)m(dw)\int_{\R^d}\frac{\overline{C}^i\overline{F}^{i}(w,y)}{|w-y|^{d+\alpha}}\overline{q}_{m-i}(t-s,y,z)\,m(dy)ds. \\
\end{eqnarray*}
First we look at the sum of the first and the second terms in the
expression of $\overline{q}_m(t,x,z)$.

There are two cases:

 Case 1. When $|x-z| \le t^{\frac{1}{\alpha}}$,
\begin{eqnarray*}
       &&C^{1}_m\int_{0}^{\frac{t}{2}}\int_{\R^d}\overline{p}(s,x,w)m(dw)\int_{\R^d}\frac{\overline{C}\overline{F}^2(w,y)}{|w-y|^{d+\alpha}}\overline{q}_{m-1}(t-s,y,z)\,m(dy)ds\\
       && +\sum_{i=2}^mC^{i}_m\int_{0}^{\frac{t}{2}}\int_{\R^d}\overline{p}(s,x,w)m(dw)\int_{\R^d}\frac{\overline{C}^i\overline{F}^{i}(w,y)}{|w-y|^{d+\alpha}}\overline{q}_{m-i}(t-s,y,z)\,m(dy)ds. \\
       &&\leq C^{1}_m\int_{0}^{\frac{t}{2}}\int_{\R^d}\overline{p}(s,x,w)m(dw)\int_{\R^d}\frac{\overline{F}^2(w,y)}{|w-y|^{d+\alpha}}\,dyds\overline{C}\overline{M}\tilde{C}_{1}(m-1)!K^{m-1}(\frac{t}{2})^{-\frac{d}{\alpha}}\\
       &&\textrm{ }\textrm{ }+\sum_{i=2}^mC^{i}_m\int_{0}^{\frac{t}{2}}\int_{\R^d}\overline{p}(s,x,w)m(dw)\int_{\R^d}\frac{\overline{F}^{2}(w,y)}{|w-y|^{d+\alpha}}\,dyds\overline{C}^2\overline{M}L^{i-2}\tilde{C}_{1}(m-i)!K^{m-i}(\frac{t}{2})^{-\frac{d}{\alpha}} \\
       &&\leq C^{1}_mC_t\overline{C}{\overline{M}}^{2}\tilde{C}_{1}(m-1)!K^{m-1}(\frac{t}{2})^{-\frac{d}{\alpha}}\\
\end{eqnarray*}
\begin{eqnarray*}
       &&\textrm{ }\textrm{ }+\sum_{i=2}^mC^{i}_mC_t\overline{C}^2{\overline{M}}^2\tilde{C}_{1}(m-i)!L^{i-2}K^{m-i}(\frac{t}{2})^{-\frac{d}{\alpha}}\\
       &&=\tilde{C}_{1}C_t\overline{C}{\overline{M}}^{2}\left(C^{1}_m(m-1)!K^{m-1}+\overline{C}\sum_{i=2}^mC^{i}_m(m-i)!L^{i-2}K^{m-i}\right)(\frac{t}{2})^{-\frac{d}{\alpha}}.\\
\end{eqnarray*}
There exits $ t_{13}>0$ with $t_{13}\leq t_{0}$  such that when
$0<t\leq t_{13}$, $C_t$ is small enough,
$$C_t\overline{C}{\overline{M}}^{2}\left(C^{1}_m(m-1)!K^{m-1}+\overline{C}\sum_{i=2}^mC^{i}_m(m-i)!L^{i-2}K^{m-i}\right)(\frac{1}{2})^{-\frac{d}{\alpha}}\leq\frac{1}{2}m!K^m.$$
Thus
\begin{eqnarray*}
       &&C^{1}_m\int_{0}^{\frac{t}{2}}\int_{\R^d}\overline{p}(s,x,w)m(dw)\int_{\R^d}\frac{\overline{C}\overline{F}^2(w,y)}{|w-y|^{d+\alpha}}\overline{q}_{m-1}(t-s,y,z)\,m(dy)ds\\
       && +\sum_{i=2}^mC^{i}_m\int_{0}^{\frac{t}{2}}\int_{\R^d}\overline{p}(s,x,w)m(dw)\int_{\R^d}\frac{\overline{C}^i\overline{F}^{i}(w,y)}{|w-y|^{d+\alpha}}\overline{q}_{m-i}(t-s,y,z)\,m(dy)ds \\
       &&\leq
       \frac{1}{2}\tilde{C}_{1}m!K^mt^{-\frac{d}{\alpha}}.\\
\end{eqnarray*}

Case 2. When $|x-z|\ge t^{\frac{1}{\alpha}}$. Let
$B_1=\{y\in\R^d|\,|y-z| \ge \frac{1}{10}|x-z|\},
B_2=\{w\in\R^d|\,|w-x| \ge  2^{-\frac{1}{2}}|x-z|\}$ and
$B_3=\{(w,y)\in \R^d \times \R^d|\,|y-z| <  \frac{1}{10}|x-z|,\,
|w-x| <  2^{-\frac{1}{2}}|x-z|\} $. On $B_3$, we have $|w-y| \ge
(1-\frac{1}{10}- 2^{-\frac{1}{2}} )|x-z|$.
\begin{eqnarray*}
       &&C^{1}_m\int_{0}^{\frac{t}{2}}\int_{\R^d}\overline{p}(s,x,w)m(dw)\int_{\R^d}\frac{\overline{C}\overline{F}^2(w,y)}{|w-y|^{d+\alpha}}\overline{q}_{m-1}(t-s,y,z)\,m(dy)ds\\
       && +\sum_{i=2}^mC^{i}_m\int_{0}^{\frac{t}{2}}\int_{\R^d}\overline{p}(s,x,w)m(dw)\int_{\R^d}\frac{\overline{C}^i\overline{F}^{i}(w,y)}{|w-y|^{d+\alpha}}\overline{q}_{m-i}(t-s,y,z)\,m(dy)ds \\
       &&\leq C^{1}_m\int_{0}^{\frac{t}{2}}\int_{\R^d}\overline{p}(s,x,w)m(dw)\int_{\R^d}\frac{\overline{C}\overline{F}^2(w,y)}{|w-y|^{d+\alpha}}\overline{q}_{m-1}(t-s,y,z)1_{B_1}(y)\,m(dy)ds\\
       &&\textrm{ }\textrm{ } +\sum_{i=2}^mC^{i}_m\int_{0}^{\frac{t}{2}}\int_{\R^d}\overline{p}(s,x,w)m(dw)\int_{\R^d}\frac{\overline{C}^i\overline{F}^{i}(w,y)}{|w-y|^{d+\alpha}}\overline{q}_{m-i}(t-s,y,z)1_{B_1}(y)\,m(dy)ds \\
       &&\textrm{ }\textrm{ } +C^{1}_m\int_{0}^{\frac{t}{2}}\int_{\R^d}\overline{p}(s,x,w)m(dw)\int_{\R^d}\frac{\overline{C}\overline{F}^2(w,y)}{|w-y|^{d+\alpha}}\overline{q}_{m-1}(t-s,y,z)1_{B_2}(w)\,m(dy)ds\\
       &&\textrm{ }\textrm{ } +\sum_{i=2}^mC^{i}_m\int_{0}^{\frac{t}{2}}\int_{\R^d}\overline{p}(s,x,w)m(dw)\int_{\R^d}\frac{\overline{C}^i\overline{F}^{i}(w,y)}{|w-y|^{d+\alpha}}\overline{q}_{m-i}(t-s,y,z)1_{B_2}(w)\,m(dy)ds \\
\end{eqnarray*}
\begin{eqnarray*}
       &&\textrm{ }\textrm{ } +C^{1}_m\int_{0}^{\frac{t}{2}}\int_{\R^d}\overline{p}(s,x,w)m(dw)\int_{\R^d}\frac{\overline{C}\overline{F}^2(w,y)}{|w-y|^{d+\alpha}}\overline{q}_{m-1}(t-s,y,z)1_{B_3}(w,y)\,m(dy)ds\\
       &&\textrm{ }\textrm{ } +\sum_{i=2}^mC^{i}_m\int_{0}^{\frac{t}{2}}\int_{\R^d}\overline{p}(s,x,w)m(dw)\int_{\R^d}\frac{\overline{C}^i\overline{F}^{i}(w,y)}{|w-y|^{d+\alpha}}\overline{q}_{m-i}(t-s,y,z)1_{B_3}(w,y)\,m(dy)ds \\
       &&\leq C^{1}_m\int_{0}^{\frac{t}{2}}\int_{\R^d}\overline{p}(s,x,w)m(dw)\int_{\R^d}\frac{\overline{F}^2(w,y)}{|w-y|^{d+\alpha}}\overline{M}\overline{C}\tilde{C}_{1}(m-1)!K^{m-1}10^{d+\alpha}\frac{(t-s)}{|x-z|^{d+\alpha}}\,dyds\\
       &&\textrm{ }\textrm{ } +\sum_{i=2}^mC^{i}_m\int_{0}^{\frac{t}{2}}\int_{\R^d}\overline{p}(s,x,w)m(dw)\int_{\R^d}\frac{\overline{F}^{2}(w,y)}{|w-y|^{d+\alpha}}\overline{M}\overline{C}^2L^{i-2}\tilde{C}_{1}(m-i)!K^{m-i}\frac{(t-s)}{|x-z|^{d+\alpha}}\,dyds \\
       &&\textrm{ }\textrm{ } +C^{1}_m\int_{0}^{\frac{t}{2}}\int_{\R^d}dw\tilde{C}_{1}{\overline{M}}^22^{\frac{1}{2}(d+\alpha)}\frac{s}{{|x-z|}^{d+\alpha}}\int_{\R^d}\frac{\overline{F}^2(w,y)}{|w-y|^{d+\alpha}}\overline{q}_{m-1}(t-s,y,z)\,dyds\overline{C}\\
       &&\textrm{ }\textrm{ } +\sum_{i=2}^mC^{i}_m\int_{0}^{\frac{t}{2}}\int_{\R^d}dw\tilde{C}_{1}{\overline{M}}^22^{\frac{1}{2}(d+\alpha)}\frac{s}{{|x-z|}^{d+\alpha}}\int_{\R^d}\frac{\overline{F}^{2}(w,y)}{|w-y|^{d+\alpha}}\overline{q}_{m-i}(t-s,y,z)\,dyds\overline{C}^2L^{i-2}\\
       &&\textrm{ }\textrm{ } +C^{1}_m\overline{C}\frac{L^2}{(1-\frac{1}{10}- 2^{-\frac{1}{2}})^{d+\alpha}}\int_{0}^{\frac{t}{2}}\int_{\R^d}\overline{p}(s,x,w)m(dw)\int_{\R^d}\frac{1}{|x-z|^{d+\alpha}}\overline{q}_{m-1}(t-s,y,z)\,m(dy)ds\\
       &&\textrm{ }\textrm{ } +\sum_{i=2}^{m-1}C^{i}_m\frac{L^i}{(1-\frac{1}{10}- 2^{-\frac{1}{2}})^{d+\alpha}}\int_{0}^{\frac{t}{2}}\int_{\R^d}\overline{p}(s,x,w)m(dw)\int_{\R^d}\frac{1}{|x-z|^{d+\alpha}}\overline{q}_{m-i}(t-s,y,z)m(dy)ds \\
       &&\textrm{ }\textrm{ }  +\frac{L^m}{(1-\frac{1}{10}- 2^{-\frac{1}{2}})^{d+\alpha}}\frac{1}{|x-z|^{d+\alpha}}\int_{0}^{\frac{t}{2}}\int_{\R^d}\overline{p}(s,x,w)m(dw)\int_{\R^d}\overline{p}(t-s,y,z)\,m(dy)ds \\
       &&\leq C^{1}_mC_t\overline{C}{\overline{M}}^{2}\tilde{C}_{1}(m-1)!K^{m-1}10^{d+\alpha}\frac{t}{|x-z|^{d+\alpha}}\\
       &&\textrm{ }\textrm{ }+\sum_{i=2}^mC^{i}_mC_t\overline{C}^2{\overline{M}}^2L^{i-2}\tilde{C}_{1}(m-i)!K^{m-i}10^{d+\alpha}\frac{t}{|x-z|^{d+\alpha}}\\
       &&\textrm{ }\textrm{ } +C^{1}_m\tilde{C_{1}}{\overline{M}}^22^{\frac{1}{2}(d+\alpha)}\frac{\tilde{C}_{1}(m-1)!K^{m-1}t}{{|x-z|}^{d+\alpha}}\int_{0}^{\frac{t}{2}}\int_{\R^d}dw\int_{\R^d}\frac{\overline{F}^2(w,y)}{|w-y|^{d+\alpha}}\frac{1}{2M_1}\overline{p}(t-s,z,y)\,dyds\overline{C}\\
       &&\textrm{ }\textrm{ } \,\,(\textrm{ by (3.3) in statement 1 for } n\leq m-1)\\
       &&\textrm{ }\textrm{ } +\sum_{i=2}^mC^{i}_m\tilde{C_{1}}{\overline{M}}^22^{\frac{1}{2}(d+\alpha)}\frac{\tilde{C}_{1}(m-i)!K^{m-i}t}{{|x-z|}^{d+\alpha}}\int_{0}^{\frac{t}{2}}\int_{\R^d}dw\int_{\R^d}\frac{\overline{F}^{2}(w,y)}{|w-y|^{d+\alpha}}\frac{1}{2M_1}\overline{p}(t-s,z,y)\,dyds\\
       &&\textrm{ }\textrm{ } \cdot \overline{C}^2L^{i-2}+C^{1}_m\frac{\overline{C}L^2}{(1-\frac{1}{10}-2^{-\frac{1}{2}})^{d+\alpha}}\frac{\frac{t}{2}}{|x-z|^{d+\alpha}}D_2\tilde{C}C_t(m-1)!K^{m-1}\\
       &&\textrm{ }\textrm{ } \,\,(\textrm{ by (3.4) in statement 1 for } n\leq m-1 \textrm{ and } \int_{\R^d}\overline{p}(s,x,w)m(dw)\leq D_2 )\\
       &&\textrm{ }\textrm{ } +\sum_{i=2}^{m-1}C^{i}_m\frac{L^i}{(1-\frac{1}{10}-2^{-\frac{1}{2}})^{d+\alpha}}\frac{\frac{t}{2}}{|x-z|^{d+\alpha}}D_2\tilde{C}C_t(m-i)!K^{m-i}\\
       &&\textrm{ }\textrm{ } +\frac{L^m}{(1-\frac{1}{10}- 2^{-\frac{1}{2}})^{d+\alpha}}\frac{t}{|x-z|^{d+\alpha}}D^2_2\\
\end{eqnarray*}
\begin{eqnarray*}
       &&\leq C^{1}_mC_t\overline{C}{\overline{M}}^{2}\tilde{C}_{1}(m-1)!K^{m-1}10^{d+\alpha}\frac{t}{|x-z|^{d+\alpha}}\\
       &&\textrm{ }\textrm{ }+\sum_{i=2}^mC^{i}_mC_t\overline{C}^2{\overline{M}}^2(m-i)!\tilde{C}_{1}L^{i-2}K^{m-i}10^{d+\alpha}\frac{t}{|x-z|^{d+\alpha}}\\
       &&\textrm{ }\textrm{ } +C^{1}_m\tilde{C_{1}}\overline{C}{\overline{M}}^22^{\frac{1}{2}(d+\alpha)}\frac{t}{{|x-z|}^{d+\alpha}}\tilde{C}_{1}(m-1)!K^{m-1}\frac{1}{2M_1}C_{t}\\
       &&\textrm{ }\textrm{ } +\sum_{i=2}^mC^{i}_m\tilde{C_{1}}\overline{C}^2{\overline{M}}^22^{\frac{1}{2}(d+\alpha)}\frac{t}{{|x-z|}^{d+\alpha}}\tilde{C}_{1}(m-i)!K^{m-i}\frac{1}{2M_1}C_{t}L^{i-2}\\
       &&\textrm{ }\textrm{ } +C^{1}_m\overline{C}\frac{L^2}{(1-\frac{1}{10}-2^{-\frac{1}{2}})^{d+\alpha}}\frac{t}{|x-z|^{d+\alpha}}D_2\tilde{C}C_t(m-1)!K^{m-1}\\
       &&\textrm{ }\textrm{ } +\sum_{i=2}^{m-1}C^{i}_m\frac{L^i}{(1-\frac{1}{10}-2^{-\frac{1}{2}})^{d+\alpha}}\frac{t}{|x-z|^{d+\alpha}}D_2\tilde{C}C_t(m-i)!K^{m-i}\\
       &&\textrm{ }\textrm{ } +\frac{L^mD^2_2}{(1-\frac{1}{10}- 2^{-\frac{1}{2}})^{d+\alpha}}\frac{t}{|x-z|^{d+\alpha}}\\
       &&\leq  (C^{1}_m\overline{C}{\overline{M}}^{2}\tilde{C}_{1}(m-1)!K^{m-1}10^{d+\alpha})C_t\frac{t}{|x-z|^{d+\alpha}}\\
       &&\textrm{ }\textrm{ }+(\sum_{i=2}^mC^{i}_m\overline{C}^2{\overline{M}}^2(m-i)!\tilde{C}_{1}L^{i-2}K^{m-i}10^{d+\alpha})C_t\frac{t}{|x-z|^{d+\alpha}}\\
       &&\textrm{ }\textrm{ } +(C^{1}_m\tilde{C_{1}}\overline{C}{\overline{M}}^22^{\frac{1}{2}(d+\alpha)}\tilde{C}_{1}(m-1)!K^{m-1}\frac{1}{2M_1})C_{t}\frac{t}{{|x-z|}^{d+\alpha}}\\
       &&\textrm{ }\textrm{ } +(\sum_{i=2}^mC^{i}_m\tilde{C_{1}}\overline{C}^2{\overline{M}}^22^{\frac{1}{2}(d+\alpha)}\tilde{C}_{1}(m-i)!K^{m-i}\frac{1}{2M_1})\frac{t}{{|x-z|}^{d+\alpha}}C_{t}L^{i-2}\\
       &&\textrm{ }\textrm{ } +(C^{1}_m\overline{C}\frac{L^2D_2}{(1-\frac{1}{10}-2^{-\frac{1}{2}})^{d+\alpha}}\tilde{C}(m-1)!K^{m-1})C_t\frac{t}{|x-z|^{d+\alpha}}\\
       &&\textrm{ }\textrm{ } +(\sum_{i=2}^{m-1}C^{i}_m\overline{C}\frac{L^iD_2}{(1-\frac{1}{10}-2^{-\frac{1}{2}})^{d+\alpha}}\tilde{C}(m-i)!K^{m-i})C_t\frac{t}{|x-z|^{d+\alpha}}\\
       &&\textrm{ }\textrm{ } +\overline{C}\frac{L^mD^2_2}{(1-\frac{1}{10}- 2^{-\frac{1}{2}})^{d+\alpha}}\frac{t}{|x-z|^{d+\alpha}}\\
\end{eqnarray*}
There exists $ t_{14}>0$ with $t_{14}\leq t_{0}$ such that when
$0<t<t_{14}$, $C_t$ is small enough and the sum of the first six
terms in the above $\leq
\frac{1}{4}\tilde{C}_1m!K^m\frac{t}{|x-z|^{d+\alpha}}$. For the
seventh term, we know that
$$ \overline{C}\frac{L^m}{(1-\frac{1}{10}-2^{-\frac{1}{2}})^{d+\alpha}}\leq\frac{1}{4}\tilde{C}_1m!K^m,\,\,\forall m>0.$$
Thus
\begin{eqnarray*}
      &&C^{1}_m\int_{0}^{\frac{t}{2}}\int_{\R^d}\overline{p}(s,x,w)m(dw)\int_{\R^d}\frac{\overline{C}\overline{F}^2(w,y)}{|w-y|^{d+\alpha}}\overline{q}_{m-1}(t-s,y,z)\,m(dy)ds\\
      && +\sum_{i=2}^mC^{i}_m\int_{0}^{\frac{t}{2}}\int_{\R^d}\overline{p}(s,x,w)m(dw)\int_{\R^d}\frac{\overline{C}^i\overline{F}^{i}(w,y)}{|w-y|^{d+\alpha}}\overline{q}_{m-i}(t-s,y,z)\,m(dy)ds \\
      &&\leq\frac{1}{2}\tilde{C}_1m!K^m\frac{t}{|x-z|^{d+\alpha}}.\\
\end{eqnarray*}
Combining case 1 and case 2, when $0<t\leq
t_{1}=\min(t_{11},t_{12},t_{13},t_{14})$, the sum of the first and
the second terms in the expression of $\overline{q}_m(t,x,z)$

\begin{eqnarray*}
       &&C^{1}_m\int_{0}^{\frac{t}{2}}\int_{\R^d}\overline{p}(s,x,w)m(dw)\int_{\R^d}\frac{\overline{C}\overline{F}^2(w,y)}{|w-y|^{d+\alpha}}\overline{q}_{m-1}(t-s,y,z)\,m(dy)ds\\
       && +\sum_{i=2}^mC^{i}_m\int_{0}^{\frac{t}{2}}\int_{\R^d}\overline{p}(s,x,w)m(dw)\int_{\R^d}\frac{\overline{C}\overline{F}^{i}(w,y)}{|w-y|^{d+\alpha}}\overline{q}_{m-i}(t-s,y,z)\,m(dy)ds \\
       &&\leq \frac{1}{2}\tilde{C}_{1}m!K^{m}t^{-\frac{d}{\alpha}}\left(1 \wedge\frac{t^{\frac{1}{\alpha}}}{|x-z|}\right)^{d+\alpha}.\\
\end{eqnarray*}
For the sum of the third and the fourth terms in the expression of
$\overline{q}_m(t,x,z)$
\begin{eqnarray*}
      &&C^{1}_m\int_{\frac{t}{2}}^{t}\int_{\R^d}\overline{p}(s,x,w)m(dw)\int_{\R^d}\frac{\overline{C}\overline{F}^2(w,y)}{|w-y|^{d+\alpha}}\overline{q}_{m-1}(t-s,y,z)\,m(dy)ds\\
      &&+\sum_{i=2}^mC^{i}_m\int_{\frac{t}{2}}^{t}\int_{\R^d}\overline{p}(s,x,w)m(dw)\int_{\R^d}\frac{\overline{C}^i\overline{F}^{i}(w,y)}{|w-y|^{d+\alpha}}\overline{q}_{m-i}(t-s,y,z)\,m(dy)ds, \\
\end{eqnarray*}
let $t-s=\tilde{s}$. This sum becomes
\begin{eqnarray*}
      &&C^{1}_m\int_{0}^{\frac{t}{2}}\int_{\R^d}\overline{p}(t-\tilde{s},x,w)m(dw)\int_{\R^d}\frac{\overline{C}\overline{F}^2(w,y)}{|w-y|^{d+\alpha}}\overline{q}_{m-1}(\tilde{s},y,z)\,m(dy)d\tilde{s}\\
      &&+\sum_{i=2}^mC^{i}_m\int_{0}^{\frac{t}{2}}\int_{\R^d}\overline{p}(t-\tilde{s},x,w)m(dw)\int_{\R^d}\frac{\overline{C}\overline{F}^{i}(w,y)}{|w-y|^{d+\alpha}}\overline{q}_{m-i}(\tilde{s},y,z)\,m(dy)d\tilde{s}. \\
\end{eqnarray*}
There are two cases:

Case a. When $|x-z| \le t^{\frac{1}{\alpha}}$,
\begin{eqnarray*}
      &&C^{1}_m\int_{0}^{\frac{t}{2}}\int_{\R^d}\overline{p}(t-\tilde{s},x,w)m(dw)\int_{\R^d}\frac{\overline{C}\overline{F}^2(w,y)}{|w-y|^{d+\alpha}}\overline{q}_{m-1}(\tilde{s},y,z)\,m(dy)d\tilde{s}\\
      &&+\sum_{i=2}^mC^{i}_m\int_{0}^{\frac{t}{2}}\int_{\R^d}\overline{p}(t-\tilde{s},x,w)m(dw)\int_{\R^d}\frac{\overline{C}^i\overline{F}^{i}(w,y)}{|w-y|^{d+\alpha}}\overline{q}_{m-i}(\tilde{s},y,z)\,m(dy)d\tilde{s} \\
      &&\leq C^{1}_m\int_{0}^{\frac{t}{2}}\int_{\R^d}\tilde{C}_{1}(\frac{t}{2})^{-\frac{d}{\alpha}}m(dw)\int_{\R^d}\frac{\overline{F}^2(w,y)}{|w-y|^{d+\alpha}}\overline{C}\tilde{C}_{1}(m-1)!K^{m-1}\frac{\overline{p}(\tilde{s},z,y)}{M_{1}}\,m(dy)d\tilde{s}\\
      &&\textrm{ }\textrm{ }+\sum_{i=2}^mC^{i}_m\int_{0}^{\frac{t}{2}}\int_{\R^d}\tilde{C}_{1}(\frac{t}{2})^{-\frac{d}{\alpha}}m(dw)\int_{\R^d}\frac{\overline{C}^2\overline{F}^i(w,y)}{|w-y|^{d+\alpha}}\tilde{C}_{1}(m-i)!L^{i-2}K^{m-i}\frac{\overline{p}(\tilde{s},z,y)}{M_{1}}\,m(dy)d\tilde{s} \\
      &&\leq C^{1}_m\tilde{C}_{1}(\frac{t}{2})^{-\frac{d}{\alpha}}\overline{C}\overline{M}^2\tilde{C}_{1}(m-1)!K^{m-1}\frac{1}{M_{1}}C_{t}
      +\sum_{i=2}^mC^{i}_m\tilde{C}_{1}(\frac{t}{2})^{-\frac{d}{\alpha}}\overline{C}^2\overline{M}^2\tilde{C}_{1}(m-i)!K^{m-i}\frac{1}{M_{1}}C_{t}\\
      &&\textrm{ }\textrm{ }+\sum_{i=2}^mC^{i}_m\tilde{C}_{1}(\frac{t}{2})^{-\frac{d}{\alpha}}\overline{C}^2\overline{M}^2\tilde{C}_{1}(m-i)!K^{m-i}\frac{1}{M_{1}}C_{t}\\
      &&=\tilde{C}_{1}\tilde{C}_1\overline{C}\overline{M}^2\left(\frac{t}{2})^{-\frac{d}{\alpha}}(C^{1}_m(m-1)!K^{m-1}\frac{1}{M_{1}}+\overline{C}\sum_{i=2}^mC^{i}_m(m-i)!K^{m-i}\frac{1}{M_{1}}\right)C_t\\
\end{eqnarray*}
There exists $ t_{21}>0$ with $t_{21}\leq t_{1}$ such that when
$0<t\leq t_{21}$, $C_t$ is small enough
$$\tilde{C}_1\overline{C}\overline{M}^2(\frac{1}{2})^{-\frac{d}{\alpha}}\left(C^{1}_m(m-1)!K^{m-1}\frac{1}{M_{1}}+\overline{C}\sum_{i=2}^mC^{i}_m(m-i)!K^{m-i}\frac{1}{M_{1}}\right)(\frac{1}{2})^{-\frac{d}{\alpha}}C_t\leq\frac{1}{2}\tilde{C}_1m!K^m.$$
Thus
\begin{eqnarray*}
      &&C^{1}_m\int_{0}^{\frac{t}{2}}\int_{\R^d}\overline{p}(t-\tilde{s},x,w)m(dw)\int_{\R^d}\frac{\overline{C}\overline{F}^2(w,y)}{|w-y|^{d+\alpha}}\overline{q}_{m-1}(\tilde{s},y,z)\,m(dy)d\tilde{s}\\
      &&+\sum_{i=2}^mC^{i}_m\int_{0}^{\frac{t}{2}}\int_{\R^d}\overline{p}(t-\tilde{s},x,w)m(dw)\int_{\R^d}\frac{\overline{C}^i\overline{F}^{i}(w,y)}{|w-y|^{d+\alpha}}\overline{q}_{m-i}(\tilde{s},y,z)\,m(dy)d\tilde{s} \\
      &&\leq \frac{1}{2}\tilde{C}_{1}m!K^mt^{-\frac{d}{\alpha}}.\\
\end{eqnarray*}

Case b. When $|x-z|\ge t^{\frac{1}{\alpha}}$. Let
$B_1=\{y\in\R^d|\,|y-z| \ge \frac{1}{10}|x-z|\},
B_2=\{w\in\R^d|\,|w-x| \ge  2^{-\frac{1}{2}}|x-z|\}$ and
$B_3=\{(w,y)\in \R^d \times \R^d|\,|y-z| <  \frac{1}{10}|x-z|,\,
|w-x| <  2^{-\frac{1}{2}}|x-z|\} $. On $B_3$,  we have $|w-y| \ge
(1-\frac{1}{10}- 2^{-\frac{1}{2}} )|x-z|$.

\begin{eqnarray*}
      &&C^{1}_m\int_{0}^{\frac{t}{2}}\int_{\R^d}\overline{p}(t-\tilde{s},x,w)m(dw)\int_{\R^d}\frac{\overline{C}\overline{F}^2(w,y)}{|w-y|^{d+\alpha}}\overline{q}_{m-1}(\tilde{s},y,z)\,m(dy)d\tilde{s}\\
      &&+\sum_{i=2}^mC^{i}_m\int_{0}^{\frac{t}{2}}\int_{\R^d}\overline{p}(t-\tilde{s},x,w)m(dw)\int_{\R^d}\frac{\overline{C}^i\overline{F}^{i}(w,y)}{|w-y|^{d+\alpha}}\overline{q}_{m-i}(\tilde{s},y,z)\,m(dy)d\tilde{s} \\
      &&\leq C^{1}_m\int_{0}^{\frac{t}{2}}\int_{\R^d}\overline{p}(t-\tilde{s},x,w)m(dw)\int_{\R^d}\frac{\overline{C}\overline{F}^2(w,y)}{|w-y|^{d+\alpha}}\overline{q}_{m-1}(\tilde{s},y,z)1_{B_1}(y)\,m(dy)d\tilde{s}\\
      &&\textrm{ }\textrm{ }+\sum_{i=2}^mC^{i}_m\int_{0}^{\frac{t}{2}}\int_{\R^d}\overline{p}(t-\tilde{s},x,w)m(dw)\int_{\R^d}\frac{\overline{C}^i\overline{F}^{i}(w,y)}{|w-y|^{d+\alpha}}\overline{q}_{m-i}(\tilde{s},y,z)1_{B_1}(y)\,m(dy)d\tilde{s} \\
      &&\textrm{ }\textrm{ }+ C^{1}_m\int_{0}^{\frac{t}{2}}\int_{\R^d}\overline{p}(t-\tilde{s},x,w)m(dw)\int_{\R^d}\frac{\overline{C}\overline{F}^2(w,y)}{|w-y|^{d+\alpha}}\overline{q}_{m-1}(\tilde{s},y,z)1_{B_2}(w)\,m(dy)d\tilde{s}\\
      &&\textrm{ }\textrm{ }+\sum_{i=2}^mC^{i}_m\int_{0}^{\frac{t}{2}}\int_{\R^d}\overline{p}(t-\tilde{s},x,w)m(dw)\int_{\R^d}\frac{\overline{C}^i\overline{F}^{i}(w,y)}{|w-y|^{d+\alpha}}\overline{q}_{m-i}(\tilde{s},y,z)1_{B_2}(w)\,m(dy)d\tilde{s} \\
      &&\textrm{ }\textrm{ }+ C^{1}_m\int_{0}^{\frac{t}{2}}\int_{\R^d}\overline{p}(t-\tilde{s},x,w)m(dw)\int_{\R^d}\frac{\overline{C}\overline{F}^2(w,y)}{|w-y|^{d+\alpha}}\overline{q}_{m-1}(\tilde{s},y,z)1_{B_3}(w,y)\,m(dy)d\tilde{s}\\
      &&\textrm{ }\textrm{ }+\sum_{i=2}^mC^{i}_m\int_{0}^{\frac{t}{2}}\int_{\R^d}\overline{p}(t-\tilde{s},x,w)m(dw)\int_{\R^d}\frac{\overline{C}^i\overline{F}^{i}(w,y)}{|w-y|^{d+\alpha}}\overline{q}_{m-i}(\tilde{s},y,z)1_{B_3}(w,y)\,m(dy)d\tilde{s} \\
      &&\leq C^{1}_m\int_{0}^{\frac{t}{2}}\int_{\R^d}\overline{p}(t-\tilde{s},x,w)m(dw)\int_{\R^d}\frac{\overline{F}^2(w,y)}{|w-y|^{d+\alpha}}\overline{C}\tilde{C}_{1}(m-1)!K^{m-1}10^{d+\alpha}\frac{\tilde{s}}{|x-z|^{d+\alpha}}\,m(dy)d\tilde{s}\\
      &&\textrm{ }\textrm{ }+\sum_{i=2}^mC^{i}_m\int_{0}^{\frac{t}{2}}\int_{\R^d}\overline{p}(t-\tilde{s},x,w)m(dw)\int_{\R^d}\frac{\overline{C}^2\overline{F}^2(w,y)}{|w-y|^{d+\alpha}}\tilde{C}_{1}(m-i)!L^{i-2}K^{m-i}10^{d+\alpha}\frac{\tilde{s}}{|x-z|^{d+\alpha}}\\
      &&\textrm{ }\textrm{ }\cdot m(dy)d\tilde{s} \\
      &&\textrm{ }\textrm{ }+ C^{1}_m\int_{0}^{\frac{t}{2}}\int_{\R^d}\,dw\tilde{C}_{1}\overline{M}^22^{\frac{1}{2}(d+\alpha)}\frac{t-\tilde{s}}{{|x-z|}^{d+\alpha}}\int_{\R^d}\frac{\overline{F}^2(w,y)}{|w-y|^{d+\alpha}}\overline{C}\tilde{C}_{1}(m-1)!K^{m-1}\frac{\overline{p}(\tilde{s},z,y)}{M_{1}}\,dyd\tilde{s}\\
      &&\textrm{ }\textrm{ }+\sum_{i=2}^mC^{i}_m\int_{0}^{\frac{t}{2}}\int_{\R^d}\,dw\tilde{C}_{1}\overline{M}^22^{\frac{1}{2}(d+\alpha)}\frac{t-\tilde{s}}{{|x-z|}^{d+\alpha}}\int_{\R^d}\frac{\overline{C}^2\overline{F}^2(w,y)}{|w-y|^{d+\alpha}}\tilde{C}_{1}L^{i-2}(m-i)!K^{m-i}\\
      &&\textrm{ }\textrm{ } \cdot \frac{\overline{p}(\tilde{s},z,y)}{M_{1}}\,dyd\tilde{s}\\
      &&\textrm{ }\textrm{ } +C^{1}_m\overline{C}\frac{L^2}{(1-\frac{1}{10}- 2^{-\frac{1}{2}})^{d+\alpha}}\int_{0}^{\frac{t}{2}}\int_{\R^d}\overline{p}(t-\tilde{s},x,w)m(dw)\int_{\R^d}\frac{1}{|x-z|^{d+\alpha}}\overline{q}_{m-1}(\tilde{s},y,z)\,m(dy)d\tilde{s}\\
      &&\textrm{ }\textrm{ } +\sum_{i=2}^{m-1}C^{i}_m\frac{L^i}{(1-\frac{1}{10}- 2^{-\frac{1}{2}})^{d+\alpha}}\int_{0}^{\frac{t}{2}}\int_{\R^d}\overline{p}(t-\tilde{s},x,w)m(dw)\int_{\R^d}\frac{1}{|x-z|^{d+\alpha}}\overline{q}_{m-i}(\tilde{s},y,z)m(dy) \\
      &&\textrm{ }\textrm{ } \cdot d\tilde{s}+\frac{L^m}{(1-\frac{1}{10}- 2^{-\frac{1}{2}})^{d+\alpha}}\frac{1}{|x-z|^{d+\alpha}}\int_{0}^{\frac{t}{2}}\int_{\R^d}\overline{p}(t-\tilde{s},x,w)m(dw)\int_{\R^d}\overline{p}(\tilde{s},y,z)\,m(dy)d\tilde{s} \\
\end{eqnarray*}
\begin{eqnarray*}
      &&\leq C^{1}_mC_{t}\overline{C}\overline{M}^2\tilde{C}_{1}(m-1)!K^{m-1}10^{d+\alpha}\frac{t}{|x-z|^{d+\alpha}}\\
      &&\textrm{ }\textrm{ }+\sum_{i=2}^mC^{i}_mC_{t}\overline{C}^2\overline{M}^2\tilde{C}_{1}(m-i)!K^{m-i}10^{d+\alpha}\frac{t}{|x-z|^{d+\alpha}}\\
      &&\textrm{ }\textrm{ } +C^{1}_m\tilde{C_{1}}\overline{C}{\overline{M}}^22^{\frac{1}{2}(d+\alpha)}\frac{t}{{|x-z|}^{d+\alpha}}\tilde{C}_{1}(m-1)!K^{m-1}\frac{1}{2M_1}C_{t}\\
      &&\textrm{ }\textrm{ } +\sum_{i=2}^mC^{i}_m\tilde{C_{1}}\overline{C}^2{\overline{M}}^22^{\frac{1}{2}(d+\alpha)}\frac{t}{{|x-z|}^{d+\alpha}}\tilde{C}_{1}(m-i)!L^{i-2}K^{m-i}\frac{1}{2M_1}C_{t}\\
      &&\textrm{ }\textrm{ } +C^{1}_m\frac{\overline{C}L^2D_2}{(1-\frac{1}{10}- 2^{-\frac{1}{2}})^{d+\alpha}}\frac{t}{|x-z|^{d+\alpha}}\tilde{C}C_t(m-1)!K^{m-1}\\
      &&\textrm{ }\textrm{ } +\sum_{i=2}^{m-1}C^{i}_m\frac{L^iD_2}{(1-\frac{1}{10}- 2^{-\frac{1}{2}})^{d+\alpha}}\frac{t}{|x-z|^{d+\alpha}}\tilde{C}C_t(m-i)!K^{m-i} \\
      &&\textrm{ }\textrm{ } \,\,(\textrm{ since statement 1 holds true for } n\leq m-1 )\\
      &&\textrm{ }\textrm{ } +\frac{L^mD^2_2}{(1-\frac{1}{10}-2^{-\frac{1}{2}})^{d+\alpha}}\frac{t}{|x-z|^{d+\alpha}}\\
      &&= (C^{1}_m\overline{C}\overline{M}^2\tilde{C}_{1}(m-1)!K^{m-1}10^{d+\alpha})C_{t}\frac{t}{|x-z|^{d+\alpha}}\\
      &&\textrm{ }\textrm{ }+(\sum_{i=2}^mC^{i}_m\overline{C}^2\overline{M}^2\tilde{C}_{1}(m-i)!K^{m-i}10^{d+\alpha})C_{t}\frac{t}{|x-z|^{d+\alpha}}\\
      &&\textrm{ }\textrm{ } +(C^{1}_m\tilde{C_{1}}\overline{C}{\overline{M}}^22^{\frac{1}{2}(d+\alpha)}\tilde{C}_{1}(m-1)!K^{m-1}\frac{1}{2M_1})C_{t}\frac{t}{{|x-z|}^{d+\alpha}}\\
      &&\textrm{ }\textrm{ } +(\sum_{i=2}^mC^{i}_m\tilde{C_{1}}\overline{C}^2{\overline{M}}^22^{\frac{1}{2}(d+\alpha)}\tilde{C}_{1}(m-i)!L^{i-2}K^{m-i}\frac{1}{2M_1})C_{t}\frac{t}{{|x-z|}^{d+\alpha}}\\
      &&\textrm{ }\textrm{ } +(C^{1}_m\overline{C}\frac{L^2D_2}{(1-\frac{1}{10}- 2^{-\frac{1}{2}})^{d+\alpha}}\tilde{C}(m-1)!K^{m-1})C_t\frac{t}{|x-z|^{d+\alpha}} \\
      &&\textrm{ }\textrm{ } +(\sum_{i=2}^{m-1}C^{i}_m\frac{L^iD_2}{(1-\frac{1}{10}- 2^{-\frac{1}{2}})^{d+\alpha}}\tilde{C}(m-i)!K^{m-i})C_t\frac{t}{|x-z|^{d+\alpha}}\\
      &&\textrm{ }\textrm{ } +\overline{C}\frac{L^mD^2_2}{(1-\frac{1}{10}- 2^{-\frac{1}{2}})^{d+\alpha}}\frac{t}{|x-z|^{d+\alpha}}.\\
\end{eqnarray*}
There exists $ t_{22}>0$ with $t_{21}\leq t_{1}$ such that when
$0<t\leq t_{22}$, $C_t$ is small enough and the sum of the first
six terms $\leq
\frac{1}{4}\tilde{C}_1m!K^m\frac{t}{|x-z|^{d+\alpha}}$. For the
seventh term, we know that
$$ \overline{C}\frac{L^mD^2_2}{(1-\frac{1}{10}-2^{-\frac{1}{2}})^{d+\alpha}}\leq\frac{1}{4}\tilde{C}_1m!K^m.$$
Thus
\begin{eqnarray*}
      &&C^{1}_m\int_{0}^{\frac{t}{2}}\int_{\R^d}\overline{p}(t-\tilde{s},x,w)m(dw)\int_{\R^d}\frac{\overline{C}\overline{F}^2(w,y)}{|w-y|^{d+\alpha}}\overline{q}_{m-1}(\tilde{s},y,z)\,m(dy)d\tilde{s}\\
      &&+\sum_{i=2}^mC^{i}_m\int_{0}^{\frac{t}{2}}\int_{\R^d}\overline{p}(t-\tilde{s},x,w)m(dw)\int_{\R^d}\frac{\overline{C}^i\overline{F}^{i}(w,y)}{|w-y|^{d+\alpha}}\overline{q}_{m-i}(\tilde{s},y,z)\,m(dy)d\tilde{s} \\
      &&\leq \frac{1}{2}\tilde{C}_{1}m!K^m\frac{t}{|x-z|^{d+\alpha}}.\\
\end{eqnarray*}
Combining case a and case b, when $0<t<t_{2}=\min(t_{21},t_{22})$,
the sum of the third and fourth terms of $\overline{q}(t,x,z)$
\begin{eqnarray*}
      &&C^{1}_m\int_{\frac{t}{2}}^{t}\int_{\R^d}\overline{p}(s,x,w)m(dw)\int_{\R^d}\frac{\overline{C}\overline{F}^2(w,y)}{|w-y|^{d+\alpha}}\overline{q}_{m-1}(t-s,y,z)\,m(dy)ds\\
      &&+\sum_{i=2}^mC^{i}_m\int_{\frac{t}{2}}^{t}\int_{\R^d}\overline{p}(s,x,w)m(dw)\int_{\R^d}\frac{\overline{C}^i\overline{F}^{i}(w,y)}{|w-y|^{d+\alpha}}\overline{q}_{m-i}(t-s,y,z)\,m(dy)ds \\
\end{eqnarray*}
$$\leq \frac{1}{2}\tilde{C}_{1}m!K^{m}t^{-\frac{d}{\alpha}}\left(1
\wedge\frac{t^{\frac{1}{\alpha}}}{|x-z|}\right)^{d+\alpha}.$$

Therefore, when $0<t\leq t_{2},$
$$ \overline{q}_m(t,x,z)\leq \tilde{C}_{1}m!K^{m}t^{-\frac{d}{\alpha}}\left(1\wedge\frac{t^{\frac{1}{\alpha}}}{|x-z|}\right)^{d+\alpha},$$
i.e. the statement $2$ holds for $n=m$. \qed

We know that for given $l>0$, $q_n(t,x,z)=q^{(l)}_n(t,x,z)$,
because of omission of the index $l$. However the expression of
$\overline{q}_n(t,x,z)$ doesn't depend on $l$. So the constant
$t_2$ in Theorem doesn't depend on $l$.

 By the above theorem, we have when $0<t\leq t_{2}$,
$$ \sum_{n=0}^{\infty}\frac{\overline{q}_n(t,x,z)}{n!}\leq \sum_{n=0}^{\infty}\tilde{C_{1}}K^nt^{-\frac{d}{\alpha}}\left(1 \wedge \frac{t^{\frac{1}{\alpha}}}{|x-z|}\right)^{d+\alpha}  =\tilde{C_{1}}\frac{1}{1-K}t^{-\frac{d}{\alpha}}\left(1 \wedge \frac{t^{\frac{1}{\alpha}}}{|x-z|}\right)^{d+\alpha}. $$
Recall that $|q_n(t,x,z)|\leq \overline{q}_n(t,x,z)$, for any
$n\geq 0$, i.e. $|q^{(l)}_n(t,x,z)|\leq \overline{q}_n(t,x,z)$,
for any $n\geq 0$, because of the omission of the index $l$. It is
clear that $\sum_{n=0}^{\infty}\frac{q^{(l)}_n(t,x,z)}{n!}$ is
uniformly convergent on $[\epsilon,t_{0}]\times\R^d\times\R^d$,
for any $\epsilon >0$. Note that $\lim_{l\rightarrow
\infty}q^{(l)}_n(t,x,z)$ exists, for any $n\geq 0$. We define
$q^{(l)}(t,x,z)=\sum_{n=0}^{\infty}\frac{q^{(l)}_n(t,x,z)}{n!}$.
Then it is easy to see that $\lim_{l\rightarrow
\infty}q^{(l)}(t,x,z)$ exists and its absolute value $\leq
\overline{q}(t,x,z)$. We denote $\lim_{l\rightarrow
\infty}q^{(l)}(t,x,z)$ by $q(t,x,z)$ and we have
\begin{equation}
       |q(t,x,z)|\leq \tilde{C_{1}}\frac{1}{1-K}t^{-\frac{d}{\alpha}}\left(1 \wedge
       \frac{t^{\frac{1}{\alpha}}}{|x-z|}\right)^{d+\alpha}.
\end{equation}

Choosing $l$ to be positive integers, we know that
$$A^{(l)}(t)=\sum_{s\leq
t}F(X_{s-},X_s)1_{\{|X_{s-}-X_{s}|>1/l\}}-\int_{0}^{t}\int_{\R^d}\frac{2C(X_s,y)F(X_s,y)1_{\{|X_s-y|>1/l\}}}{|X_s-y|^{d+\alpha}}\,m(dy)ds$$
 is a local martingale. Define
$$L^{(l)}(t)=\exp(\sum_{s\leq t}(\ln(1+F(X_{s-},X_{s})1_{\{|X_{s-},X_{s}|>1/l\}})-\int_{0}^{t}\int_{\R^d}\frac{2C(X_s,y)F(X_s,y)1_{\{|X_s-y|>1/l\}}}{|X_s-y|^{d+\alpha}}\,m(dy)ds),$$
and
$$L(t)=\exp(\lim_{l\rightarrow \infty}A^{(l)}(t)+\sum_{s\leq
      t}(\ln(1+F)-F)(X_{s-},X_{s})).$$
By an argument similar to that used in \cite{S3}, we have the
following: the limit $\lim_{l\rightarrow \infty}A^{(l)}(t)$ exists
in the sense of convergence of the norm of the space of
square-integrable martingale and convergence in probability under
$P_{x}$ and the limit is also a martingale additive functional of
$X$; $L^{(l)}(t)$ converges to $L(t)$ in probability which is a
local martingale and thus a supermartingale multiplicative
functional of $X$; $L(t)$ defines a family of probability of
measures $\{\mathbb{P}_{x}, x\in \R^d\}$ on ${\cal M}_{\infty}$ by
$d\tilde{\mathbb{P}}_{x}=L(t)d \mathbb{P}_{x}$ on ${\cal M}_{t}$.
Let $\tilde{X}=(\tilde{X}_t, \cal{M},$$\,{\cal M}_t,
\tilde{\mathbb{P}}_x, x \in {\R^d} )$ denote this new process.
Although $\tilde{X}_{t}(\omega)=X_t(\omega)$, we use $\tilde{X}_t$
to emphasize it is under $\tilde{\mathbb{P}}_{x}$. This process is
a purely discontinuous Girsanov transform of $X$.

Since $L^{(l)}(t)$ converges to $L(t)$ in probability, there
exists a subsequence $L^{(l_{i})}(t)$ of $L^{(l)}(t)$ such that
$$\lim_{l_{i}\rightarrow \infty }L^{(l_{i})}(t)=L(t).$$
On the other hand, since $\lim_{l\rightarrow
\infty}q^{(l_{i})}(t,x,z)= q(t,x,z)$ and $q(t,x,z)$ is integrable
under the measure $m(dy)$ by (3.6), we have
$$\int_{\R^d}q(t,x,z)g(z)m(dz)=\mathbb{E}_{x}[L(t)g(X_t)],$$
for any bounded measurable function $g$, i.e. $q(t,x,z)$ is the
transition density function of the purely discontinuous Girsanov
transform of $X$. Thus the upper bound of the transition density
function is obtained.

Next we establish the estimate of the lower bound of $q(t,x,z)$.
\begin{lemma}
     For any positive integer $l$, there exist $t_{2}>0$ and a large enough positive integer $k$ such that when $0<t\leq t_{2}$,
$$\frac{|q^{(l)}_{1}(t,x,z)|}{k}\leq \frac{1}{2}{p}(t,x,z),\,\, \forall (x,z) \in
\R^d\times\R^d.$$
\end{lemma}
\pf By Theorem 3.3, for any positive integer $l$, there exists
$t_{2}>0$ which doesn't depend on $l$ such that when $0<t\leq
t_{2}$,
$$|q^{(l)}_{1}(t,x,z)|\leq \overline{q}_{1}(t,x,z) \leq \tilde{C}_{1}Kt^{-\frac{d}{\alpha}}\left(1\wedge\frac{t^{\frac{1}{\alpha}}}{|x-z|}\right)^{d+\alpha}.$$
On the other hand, by (1.1) we have
$$M_{1}t^{-\frac{d}{\alpha}}\left(1\wedge\frac{t^{\frac{1}{\alpha}}}{|x-z|}\right)^{d+\alpha}\leq {p}(t,x,z).$$
Thus it is easy to see that the statement holds. \qed.

The above lemma implies that for any positive integer $l$, when
$0<t\leq t_{2}$,
$${p}(t,x,z)+\frac{q^{(l)}_{1}(t,x,z)}{k}\geq {p}(t,x,z)-\frac{|q^{(l)}_{1}(t,x,z)|}{k}\geq \frac{1}{2}{p}(t,x,z).$$
We know
$\int_{\R^d}\frac{q^{(l)}_{1}(t,x,z)}{k}\,m(dz)=\mathbb{E}_x[\frac{A^{(l)}_{t}}{k}g(X_t)],$
for any $g$ bounded measurable. Let $B_r$ be the ball $\{x\in
\R^d: \,|x|\leq r\}$. Since $1+\frac{A^{(l)}_{t}}{k} \leq
e^{\frac{A^{(l)}_{t}}{k}}$, we have
$$
 \frac{1}{|B_{r}|}\mathbb{E}_x[(1+\frac{A^{(l)}_{t}}{k})1_{B_{r}}(X_{t})] \leq \frac{1}{|B_{r}|}\mathbb{E}_x[e^{\frac{A^{(l)}_{t}}{k}}1_{B_{r}}(X_{t})]\\
$$
Thus
\begin{eqnarray*}
    \frac{1}{2}\frac{1}{|B_{r}|}\mathbb{E}_x[1_{B_{r}}(X_{t})] &\leq& \frac{1}{|B_{r}|}\mathbb{E}_x[e^{\frac{A^{(l)}_{t}}{k}}1_{B_{r}}(X_{t})]\\
                                                      &\leq&(\frac{1}{|B_{r}|}\mathbb{E}_x[e^{A^{(l)}_{t}}1_{B_{r}}(X_{t})])^{\frac{1}{k}}(\frac{1}{|B_{r}|}\mathbb{E}_x[1_{B_{r}}(X_{t})])^{1-\frac{1}{k}}\\
&&(\textrm{ by H\"older inequality }).\\
\end{eqnarray*}
Therefore
\begin{eqnarray*}
\frac{
\frac{1}{2}\frac{1}{|B_{r}|}\mathbb{E}_x[1_{B_{r}}(X_{t})]}{(\frac{1}{|B_{r}|}\mathbb{E}_x[1_{B_{r}}(X_{t})])^{1-\frac{1}{k}}}
&\leq& (\frac{1}{|B_{r}|}\mathbb{E}_x[e^{A^{(l)}_{t}}1_{B_{r}}(X_{t})])^{\frac{1}{k}}\\
\textrm{ i.e. }\\
    \frac{1}{2}(\frac{1}{|B_{r}|}\mathbb{E}_x[1_{B_{r}}(X_{t})])^{\frac{1}{k}} &\leq& (\frac{1}{|B_{r}|}\mathbb{E}_x[e^{A^{(l)}_{t}}1_{B_{r}}(X_{t})])^{\frac{1}{k}}\\
\textrm{ i.e. }\\
    \frac{1}{2^{k}}\frac{1}{|B_{r}|}\mathbb{E}_x[1_{B_{r}}(X_{t})] &\leq&
    \frac{1}{|B_{r}|}\mathbb{E}_x[e^{A^{(l)}_{t}}1_{B_{r}}(X_{t})].
\end{eqnarray*}
Let $r\downarrow 0$, we have when $0<t\leq t_{2}$
$$ \frac{1}{2^{k}}{p}(t,x,z) \leq q^{(l)}(t,x,z),\,\,\forall l>0.$$
Since $q(t,x,z)=\lim_{l\rightarrow \infty}q^{l}(t,x,z),$ we have
$$ \frac{1}{2^{k}}{p}(t,x,z) \leq q(t,x,z).$$
Thus we obtain the lower bound of $q(t,x,z)$.

It is easy to see that the semigroup property of $e^{A^{(l)}_{t}}$
implies that
$$\int_{\R^d}q^{(l)}(t,x,y)q^{l}(s,y,z)\,dy\,=q^{(l)}(t+s,x,z),\,\,\forall t,s > 0.$$
Thus by taking limit $\lim_{l\rightarrow \infty}$, we have
$$\int_{\R^d}q(t,x,y)q(s,y,z)\,dy\,=q(t+s,x,z),\,\,\forall t,s > 0.$$
By this property and the lower and upper bounds of $q(t,x,z)$ when
$0<t\leq t_{2}$, we have the following theorem
\begin{thm}
There exist positive constants $C_{3}, C_{4}, C_{5}$ and $C_{6}$
such that
\begin{equation}
C_{3}e^{-C_{4}t}t^{-\frac{d}{\alpha}}\left(1 \wedge
\frac{t^{\frac{1}{\alpha}}}{|x-z|}\right)^{d+\alpha} \leq q(t,x,z)
\leq  C_{5}e^{C_{6}t}t^{-\frac{d}{\alpha}}\left(1 \wedge
\frac{t^{\frac{1}{\alpha}}}{|x-z|}\right)^{d+\alpha}
\end{equation}
for all $(t,x,z)\in (0,\infty) \times \R^d \times \R^d $.
\end{thm}

\bigskip
\bigskip

{\bf Acknowledgement}: I am very grateful to my advisor Professor
Renming Song for his encouragement and many suggestions.

\vspace{.5in}
\begin{singlespace}

\end{singlespace}

\end{doublespace}

\end{document}